\documentclass[twoside,a4paper]{article} % double-sided
\usepackage{array}
\usepackage{booktabs}
\usepackage{amssymb}
\usepackage{a4wide} % gives a4 paper format (compared to US paper)
\usepackage{latexsym} % symbol and maths package
\usepackage{amsmath} % symbol and maths package
\usepackage{theorem} % for theorem environment in maths
\usepackage{mathrsfs} %nice calligraphic characters   %special blackboard bold
\usepackage{xspace}  % useful for macro definitions
\usepackage{fancyhdr} % for header and footer on pages
\usepackage{float} % for graphics/including pictures
\usepackage{listings}% very useful aid for including computer code in the document
\usepackage{enumerate}% useful add on for numbered lists
\usepackage{pgfplotstable}
\usepackage{morefloats}
\usepackage{tikz} 
\usepackage{accents}
\usepackage{graphicx} 
\usepackage{caption}
\usepackage{rotating}
\usepackage{subcaption}
\usepackage{algorithm}
\usepackage{algpseudocode}
\definecolor{SussexFlint}{rgb}{.00,.19,.21}
\definecolor{SussexGrey}{rgb}{.51,.58,.49}
\definecolor{SussexOrange}{rgb}{.94,.29,.00}
\definecolor{SussexYellow}{rgb}{1.00,.73,.00}
\definecolor{SussexRed}{rgb}{.94,.01,.49}
\definecolor{SussexPurple}{rgb}{.48,.06,.44}
\definecolor{SussexGreen}{rgb}{.00,.58,.46}
\definecolor{SussexBlue}{rgb}{.00,.58,.65}
\colorlet{a}{SussexOrange}
\colorlet{b}{SussexRed}
\colorlet{c}{SussexYellow}
\colorlet{d}{SussexPurple}
\colorlet{e}{SussexGreen}
\colorlet{f}{SussexBlue}
\colorlet{g}{SussexGrey}
\colorlet{h}{white}
\colorlet{i}{black}
\colorlet{j}{SussexFlint}
\numberwithin{equation}{section}
\newtheorem{theorem}{Theorem}[section]
\newtheorem{lemma}[theorem]{Lemma}
\newtheorem{corollary}[theorem]{Corollary}
\newtheorem{remark}[theorem]{Remark}
\newcommand{\qed}{\hfill$\square$}

\newcommand{\naturals}{\ensuremath{\mathbb{N}}}

\newcommand{\jump}[1]{\ensuremath{[\![#1]\!]} }

\newcommand{\Th}{\ensuremath{\mathscr{T}_h}}

\newcommand{\Ei}{\ensuremath{\mathscr{E}^i_h}}

\newcommand{\X}{\ensuremath{\Omega}}

\newcommand{\Rdsym}{\ensuremath{\mathbb R^{d\times d}_{\operatorname{Sym}}}}

\title{\begingroup\color{black}ADAPTIVE $C^0$ INTERIOR PENALTY METHODS FOR HAMILTON--JACOBI--BELLMAN EQUATIONS WITH CORDES COEFFICIENTS\endgroup}%Adaptive C
\author{Susanne C. Brenner\footnote{Susanne C. Brenner was supported in part by the National Science Foundation under Grant Nos. DMS-16-20273 and DMS-19-13035.}~ and Ellya L. Kawecki}
\begin{document}
\maketitle
\begin{abstract}
\begingroup\color{black}In this paper we conduct a priori and a posteriori error analysis of the $C^0$ interior penalty method for Hamilton--Jacobi--Bellman equations, with coefficients that satisfy the Cordes condition. These estimates show the quasi-optimality of the method, and provide one with an adaptive finite element method. In accordance with the proven regularity theory, we only assume that the solution of the Hamilton--Jacobi--Bellman equation belongs to $H^2$.\endgroup% We introduce a second method, that is adjoint consistent with a fourth-order anisotropic simply supported plate problem, proving higher rates of convergence in lower order norms, in the case that the coefficients possess higher regularity.
\end{abstract}

\begin{section}{Introduction}
The goal of this paper is to conduct a priori and a posteriori error analysis of the $C^0$ interior penalty finite element method (FEM) for the approximation of strong solutions of the following nondivergence form Hamilton--Jacobi--Bellman Dirichlet boundary-value problem. Find $u:\X\to\mathbb R$ such that
\begin{align}
\sup_{\alpha\in\Lambda}\{A^\alpha:D^2u - f^\alpha\} &= 0\quad\mbox{a.e. in }\X,\label{1}\\
u & = g\quad\mbox{on }\partial\X,\label{1a}
\end{align}
where $\X\subset\mathbb R^d$, $d\ge 2$ is convex, and $g$ is the restriction of a given $H^2(\X)$ function to $\partial\X$. We assume that 
\begin{equation}\label{0}
\Lambda\mbox{ is a compact metric space, and } A,f\in C(\overline{\X}\times\Lambda),
\end{equation} which in turn define the collection of functions $\{f^\alpha\}_{\alpha\in\Lambda},\{A^\alpha\}_{\alpha\in\Lambda}$ as follows: for each $\alpha\in\Lambda$, $f^\alpha:x\mapsto f(x,\alpha)$, $A^\alpha:x\mapsto A(x,\alpha)$. We assume that the defined collection of coefficients is uniformly elliptic in the following sense:
\begin{equation}\label{2}
\mu_1|\xi|^2\le\xi^TA^\alpha\xi\le\mu_2|\xi|^2\quad\mbox{a.e. in }\X,\,\forall\xi\in\mathbb R^d,\quad\forall\alpha\in\Lambda,
\end{equation}
as well as satisfying the following Cordes condition~\cite{MR0091400} uniformly in $\alpha$: there exists $\varepsilon\in(0,1]$ such that
\begin{equation}\label{3}
\frac{|A^\alpha|}{\operatorname{Tr}(A^\alpha)}\le\frac{1}{\sqrt{d-1+\varepsilon}}\quad\mbox{a.e. in }\X\quad\forall\alpha\in\Lambda.
\end{equation}
In the case that $\Lambda$ is a singleton set, we simply assume that $A\in L^\infty(\X)$ satisfies~(\ref{2})-(\ref{3}), and $f\in L^2(\X)$. In this case~(\ref{1})-(\ref{1a}) becomes the following linear nondivergence form elliptic equation
\begin{align}
A:D^2u & = f\quad\mbox{a.e. in }\X,\label{2a}\\
u & = g\quad\mbox{on }\partial\X.\label{2aa}
\end{align}
Remarkably, in two dimensions, uniform ellipticity implies that the the Cordes condition~(\ref{3}) (cf.~\cite{MR2260015}).%: there exists $\varepsilon\in(0,1]$ such that
%\begin{equation}\label{3}
%\frac{|A|}{\operatorname{Tr}(A)}\le\frac{1}{\sqrt{1+\varepsilon}}\quad\mbox{a.e. in }\X.
%\end{equation}

A solution $u$ of~(\ref{1})-(\ref{1a}) is called a strong solution if it belongs to $H^2(\X)$, i.e., the weak derivatives of $u$ up to second order belong to $L^2(\X)$. This means that~(\ref{1}) holds a.e. with respect to the Lebesgue measure. The linear problem~(\ref{2a})-(\ref{2aa}) is of interest, as it arises in the linearisation of~(\ref{1})-(\ref{1a}), as well as other fully nonlinear elliptic partial differential equations, such as the Monge--Amp\`{e}re (MA) equation. The MA equation, and~(\ref{1}) encompass a variety of modern applications, such as differential geometry, engineering, finance, economics, and stochastic optimal control problems.

\begingroup\color{black}\textbf{Regularity:} Since each $A^\alpha\in L^\infty(\X;\Rdsym)$, under the current hypotheses, in general a strong solution $u\in H^2(\X)$ may not belong to $H^s(\X)$ for any $s>2$. As such, we shall only assume that the true solution $u\in H^2(\X)$. 

One should note that under different hypotheses on the behaviour of the data $A$, $f$, and $\partial\X$, the solution of the linear problem~(\ref{2a})-(\ref{2aa}) may possess higher Sobolev regularity, and integrability, and may even be classically differentiable. 
\begin{itemize}
\item Calderon--Zygmund theory of strong solutions~\cite{MR1814364}: if $A\in C^0(\overline{\X};\mathbb R^{d\times d})$, %\footnote{check if symmetry is needed}
  $f\in L^p(\X)$, $1<p<\infty$, and $\partial\X\in C^{1,1}$, then $u\in W^{2,p}(\X)$.
    \item Classical solutions: if $A\in C^{0,\alpha}(\overline{\X};\mathbb R^{d\times d})$, $\alpha\in(0,1)$, $f\in C^{0,\alpha}(\overline{\X})$ and $\partial\X\in C^{2,\alpha}$, then $u\in C^{2,\alpha}(\overline{\X})$.
\end{itemize}
  %\footnote{check}, then $u\in C^{2,\alpha}(\overline{\X})$ (see~\cite{MR1814364}). We stray from such assumptions, as in general, the coefficients arising from the linearisation of nonlinear elliptic equations do not possess such continuity and regularity assumptions.
  \endgroup
 The fully nonlinear problem~(\ref{1})-(\ref{1a}) may also admit classical solutions, again provided that $A$, $f$, and $\partial\X$ are sufficiently regular. In particular, if $A,f,\partial\X\in C^\infty(\X)$ then $u\in C^0(\overline{\X})\cap C^{2,\alpha}(\X)$ for some $\alpha>0$ (cf.~\cite{MR678347}, Theorem 1, and note that $A$ is not required to satisfy~(\ref{3})). See also~\cite{MR2179357}. We seek to avoid such assumptions, as polytopal domains do not possess such regularity, and in linearising~(\ref{1})-(\ref{1a}), we cannot in general hope that the coefficients will have these properties either. See~\cite{gallistl2018numerical,MR3985933,kawecki2019curved,MR3936898} for finite element methods approximating elliptic equations on curved domains.

The main challenge in designing a numerical method for~(\ref{1})-(\ref{1a}) (aside from the nonlinearity) is the nondivergence form structure of the equation. Upon linearising~(\ref{1})-(\ref{1a}), one arrives at a sequence of problems of the form~(\ref{2a})-(\ref{2aa}). However, in general one cannot express $A:D^2u=\nabla\cdot(A\nabla u)-(\nabla\cdot A)\cdot\nabla u$, as $A\in L^\infty(\X)$, and thus may not possess sufficient regularity. This means that~(\ref{2a})-(\ref{2aa}) (and resultingly~(\ref{1})-(\ref{1a})) does not possess a weak formulation, and so, one cannot base a finite element method on that weak formulation{}{}. That said, this has not stopped the development of numerical methods for~(\ref{1})-(\ref{1a}) and~(\ref{2a})-(\ref{2aa}), often relying on the existence and uniqueness theory of the underlying equation, with methods dependent upon the different assumptions upon the coefficients and data, domain boundary, and resulting solution regularity outlined above. In particular, when
 $A\in L^\infty(\X;\mathbb R^{d\times d})$, $f\in L^2(\X)$, and $\X$ is convex, one has~\cite{gallistl2017variational,MR3985933,wang2018primal}, and if
$A\in C^0(\overline{\X};\mathbb R^{d\times d})$, $f\in L^p(\X)$, $1<p<\infty$, and $\partial\X\in C^{1,1}$ one has~\cite{feng2017finite,feng2018interior}.

The papers~\cite{MR3985933,MR3077903} present and analyse discontinuous Galerkin FEMs that utilise a discrete analogue of the Miranda--Talenti estimate; the current paper utilises a similar approach{}{}. However, the method of this paper does not involve the inclusion of additional bilinear forms which numerically enforce a discrete Miranda--Talenti estimate (as in~\cite{MR3985933,MR3077903}), and thus is simpler to implement.

The approach of~\cite{gallistl2017variational} is a mixed FEM, also relying on{}{} a variant of the Miranda--Talenti estimate, in this paper, the author was successful in proving a priori and a posteriori error estimates, as well as convergence of the adaptive method. This approach was further extended to the nonlinear setting of~(\ref{1})-(\ref{1a}) in~\cite{MR3924618}. % In this paper, we propose two finite element methods, one in the case of $L^\infty(\X)$ coefficients, and a second method that is adjoint consistent with a dual-type fourth order problem; in the latter case we are able to prove higher rates of convergence (if the coefficients are assumed to be more regular) in the $H^1$-seminorm and $L^2$-norm, similar to that of~\cite{wang2018primal}, which utilises weak discrete weak Hessian that commutes with local projections. 

The papers~\cite{feng2017finite,feng2018interior} both employ a numerical analogue of the freezing of coefficients technique utilised in the Calderon--Zygmund theory of strong solutions to~(\ref{1})-(\ref{1a}), however, the method of the present paper allows for more general coefficients and domains. For FEMs approximating~(\ref{1}) with oblique boundary conditions, see~\cite{gallistl2018numerical,MR3936898}.

The fully nonlinear setting of~(\ref{1})-(\ref{1a}) has seen several advancements in the literature, in the elliptic case~\cite{MR3623696,MR3924618,MR3307560,MR3033005,MR3196952,neilan:MT}, as well as the parabolic setting~\cite{Smears2016}. The most recent development (to the knowledge of the authors)~\cite{neilan:MT} relies on a discrete Miranda--Talenti estimate for continuous finite element functions. 

The following estimates
\begin{align}
\|D^2v\|_{L^2(\X)} &\le \|\Delta v\|_{L^2(\X)},\quad\forall v\in H^2(\X)\cap H^1_0(\X)\label{MTeq}\\
\|v\|_{H^2(\X)} &\le C \|\Delta v\|_{L^2(\X)},\forall v\in H^2(\X)\cap H^1_0(\X)\label{MTeq2}
\end{align}
are the so called Miranda--Talenti estimates, and hold when the domain $\X$ is convex. The approaches of~\cite{MR3196952,neilan:MT} rely upon renormalising the nonlinear problem with the following parameter, 
\begin{equation}\label{gammadef}
\gamma^\alpha:=\frac{A^\alpha:I}{A^\alpha:A^\alpha}\in L^\infty(\X),
\end{equation}
for each $\alpha\in\Lambda$.
%We extend this estimate (as appropriate) to functions belonging to the broken Sobolev space $H^2(\X;\Th)$ on a shape-regular triangulations $\Th$. 
%We utilise the discrete Miranda--Talenti estimate in order to prove the well-posedness of discrete schemes defined via $a_h$, on any finite dimensional subspace of $H^2(\X;\Th)$, providing one with interior penalty methods for~(\ref{1})-(\ref{1a}).

Theorem 3 of~\cite{MR3196952} provides the existence and uniqueness of a function $u$ belonging to the space
$$H:=H^2(\X)\cap H^1_0(\X),$$ that satisfies~(\ref{1})-(\ref{1a}), in the case that $g\equiv 0$. 
\begingroup\color{black}Treating the case of inhomogeneous boundary data follows in a manner similar to that of~\cite{MR3196952}, Theorem 3. With the aim of invoking the Browder--Minty Theorem, we first define $F_\gamma:H^2(\X)\to L^2(\X)$ by 
\begin{equation}\label{Fgammadef}
F_\gamma[u]:=\sup_{\alpha\in\Lambda}\{\gamma^\alpha(A^\alpha:D^2u-f^\alpha)\},
\end{equation}
and proceed to define $a_g:H\to H'$ (where $H'$ denotes the dual space of $H$) by
\begin{equation}\label{a:def}
a_g(u;v):=(F_\gamma[u+g],\Delta v)_{L^2(\X)}\quad u,v\in H.
\end{equation}
One can show that $a_g$ is strictly monotone, and Lipschitz continuous on $H$, yielding the existence and uniqueness of a function $u_0\in H$ such that 
\begin{equation}a_g(u_0;v)=0\quad\forall v\in H. \label{variation}
\end{equation}
Finally, we uniquely define $u:=u_0+g$, which satisfies~(\ref{1})-(\ref{1a}). This provides us with the following theorem.  \begin{theorem}\label{inhomog:eandu:15:11}
Assume that $\X\subset\mathbb R^d$ is a convex domain, and that the collection $\{A^\alpha\}_{\alpha\in\Lambda}$ satisfies~(\ref{2})-(\ref{3}). Furthermore, assume that $g\in H^{2}(\X)$. Then, there exists a unique strong solution $u\in H^2(\X)$ of the following HJB equation:
\begin{equation}\label{inhomog:cordes:HJB}
\begin{aligned}
\sup_{\alpha\in\Lambda}\{A^\alpha:D^2 u-f^\alpha\} & = 0\quad\mbox{a.e. in}\,\,\X,\\
u & = g\quad\mbox{on}\,\,\partial\X.
\end{aligned}
\end{equation}
\end{theorem}
\endgroup

\textbf{Contributions:} In this paper we obtain a priori and a posteriori error estimates under the assumption that the true solution belongs to $H^2(\X)$. We note that the method we present has been considered in~\cite{neilan:MT}, in the homogeneous Dirichlet case, where the author proves stability, and a priori error estimates for the problem~(\ref{1})-(\ref{1a}), as well as the fully nonlinear Hamilton--Jacobi--Bellman equation. Our approach to the stability analysis is distinct from that of~\cite{neilan:MT}, as we also consider the case of inhomogeneous boundary data. Furthermore, the recent publication~\cite{MR4019733} provides the existence of an enrichment operator when $p\ge 2$, and $d\in\{2,3\}$, which  removes the restriction upon the polynomial degree $p\in\{2,3\}$, when $d=3$ present in~\cite{neilan:MT} (cf.~\cite{neilan:MT} Remark 4). Furthermore, we also undertake a posteriori error analysis for this problem, and justify that one may utilise the scheme to approximate solutions to the fully nonlinear Monge--Amp\`{e}re equation (see Section~\ref{MA:apps}).

As mentioned, a motivation of this paper is to develop a numerical method for the Monge--Amp\`{e}re (MA) equation. The (MA) equation is a prototypical fully nonlinear elliptic equation, arising in differential geometry, optimal transport, engineering and fluid dynamics: given $f:\X\to\mathbb R^+$, uniformly positive, and $g:\partial\X\to\mathbb R$, find $u:\X\to\mathbb R$ such that
\begin{equation}\label{MAD:star}
\left\{
\begin{aligned}
\operatorname{det}D^2u & = f\quad\mbox{in }\X,\\
u & = g\quad\mbox{on }\partial\X.
\end{aligned}
\right.
\end{equation}
%However, we observe that in the two dimensional case, with homogeneous boundary conditions ($g\equiv0$),  a convex function $u\in W^{2,\infty}(\X)$ satisfies
%\begin{equation}\label{2dexample16:09}
%\left\{
%\begin{aligned}
%\operatorname{det} D^2u(x)& = f(x),\,\,x\in\X,\\
%u(x) & = 0,\,\,x\in\partial\X,
%\end{aligned}
%\right.
%\end{equation}
%then one can see that $\tilde{u}:=-u$ is \emph{concave}, and satisfies
%\begin{equation*}
%\left\{
%\begin{aligned}
%\operatorname{det} D^2\tilde u(x)& = \operatorname{det}(-D^2u(x)) = \operatorname{det}(D^2u(x)) =  f(x),\,\,x\in\X,\\
%\tilde u(x) & = -u(x) = 0,\,\,x\in\partial\X,
%\end{aligned}
%\right.
%\end{equation*}
%i.e., $\tilde{u}$ also solves~(\ref{2dexample16:09}). 
In general~(\ref{MAD:star}) may admit at most two solutions; a simple example is given in the case $d=2$, $g\equiv0$, where it is clear that if $u$ satisfies~(\ref{MAD:star}), then so does $-u$.
One would hope that numerical methods that approximate solutions of~(\ref{MAD:star}) may have the same uniqueness property, that is, that there exists at most two solutions to the numerical method. However, this is not always the case. In~\cite{MR3049920}, the authors implement a standard nine-point stencil finite difference method for the problem~(\ref{MAD:star}) with a smooth right-hand side, and a smooth solution $u$, with the choice of domain $\X=(0,1)^2$. Upon implementing this method on a $4\times 4$ grid, solving the resulting nonlinear system by applying Newton's method, they obtain sixteen different numerical solutions by varying the initial guess of the Newton's method. 

As mentioned in~\cite{MR3049920}, one may conjecture that this phenomena extrapolates, causing Newton's method to potentially converge to $2^{(N-2)^2}$ different solutions on and $N\times N$ grid, by varying the initial guess. When designing a numerical scheme, it is important that one knows which solution the method is converging to, without needing too much (Newton's method is well known to be conditionally convergent, and a prerequisite to convergence is often sufficient proximity to the true solution) prior knowledge of the true solution. Indeed, the aforementioned finite difference method implemented in~\cite{MR3049920}, was proposed in~\cite{MR2683581}, with an additional \emph{selection} criteria, which in essence singles out a particular numerical solution.

We overcome this difficulty, by utilising a long standing result due to N. Krylov~\cite{MR901759}, which allows one to characterise the MA equation~(\ref{MAD:star}) as a HJB equation, %with the Krylov control set, $X$, defined as follows
%\begin{equation}\label{KrylCont}
%X:=\{W\in\mathbb R^{d\times d}_{\operatorname{Sym}}:W\ge 0,\operatorname{Tr}W=1\}.
%\end{equation} 
%In particular,~(\ref{MAD:star}) is equivalent to the following problem: 
%\begin{equation}\label{HJB-MA-degenerate}
%\left\{
%\begin{aligned}
%\max_{W\in X}\{-W:D^2u+d(\operatorname{det} W)^{1/d}f^{1/d}\} & = 0\quad\mbox{in }\X,\\
%u & = g\quad\mbox{on }\partial\X,
%\right.
%\end{equation}
if and only if $u$ is convex. In the case that $u\in W^{2,\infty}(\X)$ is \emph{uniformly} convex, and $d=2$, we are able to further show that the resulting HJB equation is equivalent to one with a control set $\Lambda$, and data $A,f$ that satisfy~(\ref{0})-(\ref{3}). Moreover, the resulting numerical scheme is uniquely solvable. For other numerical method for the approximation of solutions to the MA problem, see~\cite{MR3020056,MR3623696,kawecki:lakkis:pryer,MR3162358}.

%The main contributions of this paper are the proposition of a $C^0$-interior penalty method for~(\ref{1})-(\ref{1a}), with stability, a priori, and a posteriori error analysis in a $H^2$-type norm. %Furthermore, we propose a second method, for which we are able to undertake the same analysis, with the addition of error estimates in the  $H^1$-seminorm and $L^2$-norm. A final novel contribution, is that in designing the second method, that is adjoint consistent, we are provided with a finite element method for the approximation solutions to the adjoint problem; a particular example of this adjoint problem is the simply supported plate problem (for numerical methods for more general versions of this problem, see~\cite{semper1992conforming} and~\cite{brenner2011c} for conforming and nonconforming methods, respectively).

\textbf{Domain assumptions:} In the current section, we have assumed that $\X\subset\mathbb R^d$, $d\ge2$, is convex, as this is a sufficient assumption of Theorem~\ref{inhomog:eandu:15:11}. However, in Section~\ref{sec:2} we provide the numerical scheme, and from this point on, we further assume that $d\in\{2,3\}$ and that $\X$ is polytopal.

\begingroup\color{black}This paper is laid out as follows. %In Section~\ref{newsec:2} we discuss the PDE analysis framework for the well posedness of~(\ref{1})-(\ref{1a}). 
In Section~\ref{sec:2}, we introduce the discrete problem, and prove the stability of the associated bilinear form. Section~\ref{sec:4} is devoted to convergence analysis; we prove quasi-optimal apriori error estimates and a posteriori error estimates in a $H^2$-type norm. In Section~\ref{iterative:sec}, we propose the linearisation scheme and adaptive scheme. Section~\ref{MA:apps} is devoted to applications to the Monge--Amp\`{e}re problem. In Section~\ref{sec:5} we implement the proposed finite element method (as well as the adaptive version) in FEniCS~\cite{AlnaesBlechta2015a}, confirming the theoretical results of the paper. Finally, in Section~\ref{sec:6}, we provide concluding remarks on what has been achieved in this paper.\endgroup
\end{section}
\begin{section}{The Discrete Problem}\label{sec:2}
As mentioned in the introduction, from this point on, we shall further assume that $\X\subset\mathbb R^d$, $d\in\{2,3\}$ is convex and polytopal.
Let $\Th$ be a simplicial triangulation of $\X$ and $V_h\subset H^1(\X)$ be the continuous Lagrange finite element space of order $p\ge2$ associated with $\Th$, and denote $V_{h,0}:=V_h\cap H^1_0(\X)$. %The proceeding results require the consideration of the discontinuous Galerkin finite element space of degree $p$, denoted by $\mathcal{V}_h$. We have that $V_h\subset\mathcal{V}_h\subset H^2(\X;\Th)$. 
We denote by $D^2_h$ and $\Delta_h$, the piecewise Hessian and Laplacian, respectively. 
Furthermore, we shall make use of the following mesh dependent (semi)norm for $u\in H^2(\X;\Th):=\{v\in L^2(\X):v|_K\in H^2(K)\,\forall K\in\Th\}$
\begin{align}
%(u,v)_h&:=\int_\X D^2_hu:D^2_hv+\sum_{e\in\Ei}\int_e\frac{\sigma}{h_e}\jump{\nabla u}_e\cdot\jump{\nabla v}_e+\sum_{e\in\Eib}\int_e\frac{\rho}{h_e^3}\jump{u}_e\jump{v}_e,\\
%(u,v)_{J_h^{\sigma,\rho}}&:=\sum_{e\in\Ei}\int_e\frac{\sigma}{h_e}\jump{\nabla u}_e\cdot\jump{\nabla v}_e+\frac{\sigma}{h_e^3}\jump{u}_e\jump{v}_e+\sum_{e\in\Eb}\int_e\frac{\rho}{h_e^3}\jump{u}_e\jump{v}_e,\\
\|u\|_h^2&:=\int_\X |D^2_hu|^2+\sum_{e\in\Ei}\frac{\sigma}{h_e}\|\jump{\partial u/\partial n}_e\|_{L^2(e)}^2,
%\|u\|_h^2&:=(u,u)_h,%(u,u)_h,\quad|u|_{J_h^{\sigma,\rho}}^2:=(u,u)_{J_h^{\sigma,\rho}}\label{star15:28},
\end{align}
and we note that $\|\cdot\|_h$ is indeed a norm on $V_{h,0}$.

The discrete problem is posed as follows: we seek  $u_h\in V_h$ satisfying
\begin{equation}\label{thescheme2}
\begin{aligned}
a_h(u_h;v) &:= \int_\X F_\gamma[u_h]\,\Delta_h v +\sum_{e\in\Ei}\frac{\sigma}{h_e}\int_e\jump{\partial u_h/\partial n}_e\jump{\partial v/\partial n}_e = 0\quad\forall v\in V_{h,0},\\
u_h|_{\partial\X}&:=g_h,
\end{aligned}
\end{equation}
where $g_h\in V_h$ is a suitable approximation of $g$ (the derivatives in $F_\gamma$ defined by~(\ref{Fgammadef}){}{} are considered piecewise), $\Ei$ is the set of internal edges of $\Th$, $\jump{\cdot}_e$ denotes the jump across an edge $e$, and $\sigma$ is a positive constant.
%where $g_h\in V_h$ is the $(\cdot,\cdot)_h$-projection of $g$ onto $V_h$ (the derivatives in $F_\gamma$ are considered piecewise). In the second case, we seek $u_h\in V$ such that
%\begin{equation}\label{thescheme}
%a_h(u_h;v)  = (g_h,v)_{J_h^{0,\rho}}\quad\forall v\in V.
%\end{equation}
%The two approaches are similar, but lend themselves naturally to different considerations of spaces $V$. In particular, when $V=V_h$, one would often consider~(\ref{thescheme2}), and when $V=\mathcal{V}_h$, one would usually employ~(\ref{thescheme}). We shall focus on~(\ref{thescheme}), and later remark on how one may prove existence and uniqueness of a solution to~(\ref{thescheme2}) (i.e., when the solution and test spaces may differ).

We remark that if $g\equiv0$, and we instead seek $u_h\in V_{h,0}$, then~(\ref{thescheme2}) coincides with the method presented in~\cite{neilan:MT}. We also note that the scheme is consistent in the following sense: if $u\in H^2(\X)$ satisfies~(\ref{1})-(\ref{1a}), then
\begin{equation}\label{cons:13:45}
a_h(u;v) = 0\quad\forall v\in V_{h,0}.
\end{equation}
The above holds, since $u$ satisfies~(\ref{1})-(\ref{1a}), and $u\in H^2(\X)$ so $\jump{\partial u/\partial n}_e=0$ for $e\in\Ei$.
The following theorem and corollary are from~\cite{neilan:MT}. As mentioned in the introduction, the results that follow, as presented in~\cite{neilan:MT} hold for $d=2$, for any $p\ge 2$, and for $d=3$ if $p\in\{2,3\}$. However, this is occurs as the proofs rely on the existence of an operator $E_h:V_{h,0}\to H^2(\X)\cap H^1_0(\X)$ (called an enrichment operator), that in particular satisfies the following estimate:
\begin{equation}\label{14:18}
\|E_hv-v\|_h\le C_*\sum_{e\in\Ei}\frac{1}{h_e}\|\jump{\partial v/\partial n}_e\|_{L^2(e)}^2,\quad\forall v\in V_{h,0}
\end{equation}
where the constant $C_*$ is (in principle) a computable, positive constant dependant only on the shape regularity of $\Th$. A particular construction of such an operator is provided in~\cite{neilan:MT} and uses the $C^1$ family of Clough--Tocher spaces, which leads to the aforementioned restriction when $d=3$ (cf.~\cite{neilan:MT}, Remark 4). However, in the recent paper~\cite{MR4019733}, the existence of an operator that satisfies~(\ref{14:18}) has been proven, only assuming $p\ge2$, for $d\in\{2,3\}$. Thus, the proceeding results hold for $d\in\{2,3\}$ and $p\ge 2$.
\begin{theorem}\label{thm:A}
One has that for any $v_h\in V_{h,0}$,
\begin{equation}\label{discmt1}
\|D^2_hv_h\|_{L^2(\X)}\le\|\Delta_h v_h\|_{L^2(\X)}+C\left(\sum_{e\in\Ei}\frac{1}{h_e}\|\jump{\partial v_h/\partial n}\|_{L^2(e)}^2\right)^{1/2},
\end{equation}
where the constant $C$ is independent of $h$. 
\end{theorem}
\begin{corollary}\label{corr:A}
One has that for any $v\in V_{h,0}$, and all $t\in(0,1)$,
\begin{equation}\label{discMT3}
\|\Delta_hv\|_{L^2(\X)}^2\ge(1-t)\|D^2_hv_h\|_{L^2(\X)}^2-\frac{C^2}{t}\left(\sum_{e\in\Ei}\frac{1}{h_e}\|\jump{\partial v/\partial n}\|_{L^2(e)}^2\right),
\end{equation}
where the constant $C$ is the constant present in Theorem~\ref{thm:A}.
\end{corollary}
%\emph{Proof:} Utilising the Cauchy--Schwarz inequality in Theorem~\ref{thm:B}, we obtain the following for any $\delta>0$,
%$$\|D^2_hv_h\|_{L^2(\X)}^2\le(1+\delta)\|\Delta_hv\|_{L^2(\X)}^2+C^2(1+1/\delta)\left(\sum_{e\in\Eione}\frac{1}{h_e}\|\jump{\nabla v}\|_{L^2(e)}^2+\sum_{e\in\Eibone}\frac{1}{h_e^3}\|\jump{v}\|_{L^2(e)}^2\right).$$
%Setting $t=\delta/(1+\delta)\in(0,1)$, and rearranging, we obtain
%\begin{align*}
%\|\Delta_hv\|_{L^2(\X)}^2 & \ge(1-t)\|D^2_hv_h\|_{L^2(\X)}^2-C^2\frac{1-t}{t}\left(\sum_{e\in\Eione}\frac{1}{h_e}\|\jump{\nabla v}\|_{L^2(e)}^2+\sum_{e\in\Eibone}\frac{1}{h_e^3}\|\jump{v}\|_{L^2(e)}^2\right)\\
%&\ge(1-t)\|D^2_hv_h\|_{L^2(\X)}^2-\frac{C^2}{t}\left(\sum_{e\in\Eione}\frac{1}{h_e}\|\jump{\nabla v}\|_{L^2(e)}^2+\sum_{e\in\Eibone}\frac{1}{h_e^3}\|\jump{v}\|_{L^2(e)}^2\right).
%\end{align*}
%\hfill$\square$\\
We now prove a strict monotonicity result for $a_h$ (a variant of~\cite{neilan:MT}, Lemma 7), provided that $\sigma$ is sufficiently large.
\begin{lemma}\label{lemma:C}
One has that for any $u,v\in V_{h}$, such that $u-v\in V_{h,0}$
$$a_h(u;u-v)-a_h(v;u-v)\ge \delta(1-\sqrt{1-\varepsilon})\|u-v\|_h^2,$$
for any $\delta\in(0,1)$, independent of $h$, $u$ and $v$, provided $\sigma$ is sufficiently large (dependent on $\delta$).
\end{lemma}
\emph{Proof:} Take $u,v$, as in the hypotheses of the lemma, and denote $w:=u-v\in V_{h,0}$. \begingroup\color{black}Denoting $I$ to be the $d\times d$ identity matrix, by~(\ref{3}), we have that\endgroup
\begin{equation}\label{gamma:2}
\begin{aligned}
|\gamma^\alpha A^\alpha-I|^2 &= (\gamma^\alpha A^\alpha-I):(\gamma^\alpha A^\alpha-I) = (\gamma^\alpha)^2(A^\alpha:A^\alpha)-2\gamma^\alpha(A^\alpha:I)+I:I\\
&=-\frac{(A^\alpha:I)^2}{(A^\alpha:A^\alpha)}+d\\
&\le-(d-1+\varepsilon)+d = 1-\varepsilon\quad\mbox{a.e. in }\X,\quad\forall\alpha\in\Lambda.
\end{aligned}
\end{equation}
Inequality~(\ref{gamma:2}), Theorem~\ref{thm:A}, and Corollary~\ref{corr:A} imply that for any $t\in(0,1)$
\begin{align*}
\int_\X(F_\gamma[u]-F_\gamma[v])\,\Delta_h w
&\ge \|\Delta_hw\|_{L^2(\X)}^2-\int_\X\sup_{\alpha\in\Lambda}\{|(\gamma^\alpha A^\alpha-I_d):D^2_hw|\}|\Delta_h w|\\
&\ge\|\Delta_hw\|_{L^2(\X)}^2-\sqrt{1-\varepsilon}\|D^2_hw\|_{L^2(\X)}\|\Delta_hw\|_{L^2(\X)}\\
&\ge(1-\sqrt{1-\varepsilon}/2)\|\Delta_hw\|_{L^2(\X)}^2-(\sqrt{1-\varepsilon}/2)\|D^2_hw\|_{L^2(\X)}^2\\
&\ge\left((1-t)(1-\sqrt{1-\varepsilon}/2)-\sqrt{1-\varepsilon}/2\right)\|D^2_hw\|_{L^2(\X)}^2-\frac{C}{t}|w|_{J_h}^2.
\end{align*}
Now, for a given $\delta\in(0,1)$, we set $t=(1-\delta)(1-\sqrt{1-\varepsilon})/(1-\sqrt{1-\varepsilon}/2)$. This gives us
\begin{align*}
a_h(u;u-v)-a_h(v;u-v)&\ge\delta(1-\sqrt{1-\varepsilon})\|w\|_{h}^2\\
&~~+\left(\sigma(1-\delta(1-\sqrt{1-\varepsilon}))-C\frac{(1-\sqrt{1-\varepsilon}/2)^2}{(1-\delta)(1-\sqrt{1-\varepsilon})}\right)\sum_{e\in\Ei}\frac{1}{h_e}\|\jump{\partial w/\partial n}\|_{L^2(e)}^2\\
&\ge\delta(1-\sqrt{1-\varepsilon})\|w\|_{h}^2
\end{align*}
provided that $\sigma$ satisfies
$$\sigma\ge \frac{C(1-\sqrt{1-\varepsilon}/2)^2}{(1-\delta)(1-\sqrt{1-\varepsilon})(1-\delta(1-\sqrt{1-\varepsilon}))}.$$
\hfill$\square$\\
The following lemma is a simple consequence of the definition of the norm $\|\cdot\|_h$.
\begin{lemma}\label{lemma:D}
%Let $\X\subset\mathbb R^d$, be polytopal. Then, for $h>0$, 
One has that for any $u,v,w\in V_{h}$,
\begin{align*}
|a_h(u;w)-a_h(v;w)|&\le C\|u-v\|_h\|w\|_h,
%\sum_{e\in\Ei}\frac{\sigma}{h_e}\int_e\jump{\partial g_h/\partial n}_e\jump{\partial v/\partial n}_e&\le\|g\|_h\|v\|_h,
\end{align*}
where the constant $C$ is independent of $u,v,w$, and $h$.
\end{lemma}
%Lemma~\ref{lemma:C} and Lemma~\ref{lemma:D} combined with the Browder--Minty Theorem imply the following Theorem.
%\begin{theorem}
%There exists a unique $u_h\in V_{h,0}$ such that
%$$a_h(u_h;v)=0\quad\forall v\in V_{h,0}.$$
%\end{theorem}
%We now prove that~(\ref{thescheme2}) is also uniquely solvable.
The following proof is motivated by the proof of Theorem~\ref{inhomog:eandu:15:11}.
\begin{theorem}
Under the hypotheses of Lemma~\ref{lemma:C}, there exists a unique $u_h\in V_{h}$ satisfying~(\ref{thescheme2}).
\end{theorem}
\emph{Proof:} Let us define $a_{g_h}:V_{h,0}\times V_{h,0}\to\mathbb R$, by $a_{g_h}(u_h;v):=a_{h}(u+g_h;v)$ for all $u_h,v\in V_{h,0}$. Lemmas~\ref{lemma:C} and~\ref{lemma:D} then imply that for all $u_h,v,w\in V_{h,0}$, and for any $\delta\in(0,1)$ (so long as $\sigma$ is sufficiently large, dependent on $\delta$)
\begin{align*}
a_{g_h}(u_h;u_h-v)-a_{g_h}(v;u_h-v)&\ge\delta(1-\sqrt{1-\varepsilon})\|u_h-v\|_h^2,\\
|a_{g_h}(u_h;w)-a_{g_h}(v;w)|&\le C\|u_h-v\|_h\|w\|_h,
\end{align*}
where the constant $C$ is independent of $u_h,v,w$. Thus $a_{g_h}$ is strongly monotone and Lipschitz continuous, and so the Browder--Minty Theorem implies the existence and uniqueness of $u_{h,0}\in V_{h,0}$ such that
$$a_{g_h}(u_{h,0},v)=a_{h}(u_{h,0}+g_h;v)=0\quad\forall v\in V_{h,0}.$$
Thus, we may uniquely define $u_h:=u_{h,0}+g_h$, which satisfies~(\ref{thescheme2}).
\hfill$\square$
\end{section}
\begin{section}{Convergence Analysis}\label{sec:4}
\subsection{A Priori Error Analysis}
For the remainder of the paper, we assume that the parameter $\sigma$ is chosen such that that exists a unique $u_h\in V_h$ satisfying~(\ref{thescheme2}). %We shall also assume that $g_h$ is the $(\cdot,\cdot)_h$-best approximation of $g$ in $V_h$.
\begingroup\color{black}\begin{remark}[Choice of $g_h$]
In practice, one may use a variety of numerical approximations $g_h$ of $g$, for example, the $L^2$ projection, or some suitable interpolant. However, for the density argument of Remark~(\ref{densityremark}){}{} it is useful to define $g_h$ to be the unique element of $V_h$ that satisfies
\begin{equation}\label{ghdef}\int_\X D^2_hg_h:D^2_hv+\int_{\partial\X}g_hv+\sum_{e\in\Ei}\frac{\sigma}{h_e}\int_e\jump{\partial g_h/\partial n}_e\jump{\partial v/\partial n}_e=\int_\X D^2_hg:D^2_hv+\int_{\partial\X}gv\quad\forall v\in V_h.
\end{equation}
\end{remark}\endgroup
We first prove a quasi-optimal error estimate for the error $\|u-u_h\|_h$, where $u\in H^2(\X)$ satisfies~(\ref{1})-(\ref{1a}). Let $v\in V_h$ satisfy $v|_{\partial\X} = g_h$, where $g_h$ satisfies~(\ref{ghdef}). The triangle inequality gives us
\begin{equation}\label{4.1}
\|u-u_h\|_h\le\|u-v\|_h+\|v-u_h\|_h.%\le\|u-v\|_h,
%+\frac{1}{\delta}\sup_{w\in V_h}\frac{a_h(v-u_h,w)}{\|w\|_h}.
\end{equation}
Lemma~\ref{lemma:C},~(\ref{cons:13:45}), and Lemma~\ref{lemma:D} imply that for any $\delta\in(0,1)$ (denoting $c_{\delta,\varepsilon}:=\delta^{-1}(1-\sqrt{1-\varepsilon})^{-1}$)
\begin{align*}
c_{\delta,\varepsilon}^{-1}\|v-u_h\|_h^2&\le a_h(u_h;u_h-v)-a_h(v;u_h-v)\\
&=a_h(u;u_h-v)-a_h(v;u_h-v)\\
&\le C\|u-v\|_h\|u_h-v\|_h.
\end{align*}
%If $u_h\in V_h$ satisfies~(\ref{thescheme}), then 
%$$a_h(u_h;u_h-v)-(g,u_h-v)_{J_h^{0,\rho}}=(g_h-g,u_h-v)_{J_h^{0,\rho}}\le|g_h-g|_{J_h^{0,\rho}}\|u_h-v\|_h,$$
%and so
%\begin{equation}\label{4.2}
%\|v-u_h\|_h\le C(\|u-v\|_h+|g_h-g|_{J_h^{0,\rho}}).
%\end{equation}
%If $u_h\in V$ satisfies~(\ref{thescheme2}), we further assert that $v|_e=g_h$ for all $e\in\Eb$, and so 
%$$a_h(u_h;u_h-v)-(g,u_h-v)_{J_h^{0,\rho}}=0,$$
Thus,
\begin{equation}\label{4.2b}
\|v-u_h\|_h\le C\|u-v\|_{h}.
\end{equation}
Combining~(\ref{4.1}) with~(\ref{4.2b}), we arrive at the following quasi-optimal error estimate.
\begin{theorem}\label{16:27}
There exists a positive constant $C_\sharp$, independent of $h$, such that if $u_h\in V_h$ satisfies~(\ref{thescheme2}), then
\begin{equation}\label{4.3}
\|u-u_h\|_h\le C_\sharp(\inf_{v\in V_h:v|_{\partial\X}=g_h}\|u-v\|_h).
\end{equation}
%Furthermore, if $u_h\in V$ satisfies~(\ref{thescheme2}), then
%\begin{equation}\label{4.3b}
%\|u-u_h\|_h\le C_\sharp(\inf_{v\in V:v|_e=g_h\,\forall e\in\Eb}\|u-v\|_h).
%\end{equation}
\end{theorem}
%\begin{corollary}
%There exists a positive constant $C_\flat$, independent of $h$, such that
%\begin{equation}
%|u-u_h|_{H^2(\X;\Th)}\le C_\flat\inf_{v\in V_h}|u-v|_{H^2(\X;\Th)}.
%\end{equation}
%\end{corollary}
%\emph{Proof:} Due to Theorem~\ref{16:27}, we have that
%\begin{equation}\label{4.3}
%|u-u_h|_{H^2(\X;\Th)}\le\|u-u_h\|_h\le C_\sharp\inf_{v\in V_h}\|u-v\|_h\le C_\sharp\left(\inf_{v\in V_h}|u-v|_{H^2(\X;\Th)}+\left(\sum_{e\in\Ei}\int_e\jump{\partial u_h/\partial n}^2\,ds\right)^{1/2}\right).
%\end{equation}
%
%Since $u\in H^2(\X)\subset\cap H^2(\X;\mathscr{T}_{h_0})$, Lemma~\ref{lemma:A} and~(\ref{4.3}) and imply the following corollary.
\begingroup\color{black}\begin{remark}\label{densityremark}
Estimate~(\ref{4.3}) in combination with a density argument shows that $$\lim_{h\to0}\|u-u_h\|_h=0.$$
Moreover, the Poincar\'e--Friedrichs inequality for piecewise $H^2$ functions (cf.~\cite{MR3606463,MR2106270}), implies that exists a positive constant $C$ independent of $h$ such that
\begin{equation}\label{4.4}
\|u-u_h\|_{H^1(\X)}+\|u-u_h\|_{L^\infty(\X)}\le C\|u-u_h\|_h,
\end{equation}
and so
$$\lim_{h\to0}(\|u-u_h\|_{H^1(\X)}+\|u-u_h\|_{L^\infty(\X)}) = 0.$$
\end{remark}\endgroup
\subsection{A Posteriori Error Analysis}
The \emph{a posteriori} error analysis is based on an enriching operator $E_h:V_{h,0}\to H^2(\X)\cap H^1_0(\X)$ that satisfies~(\ref{14:18}).
%\begin{equation}\label{4.5}
%\|v-E_hv\|_h^2\le C_*\sum_{e\in\Ei}\frac{1}{h_e}\int_e\jump{\partial v/\partial n}_e^2ds\quad\forall v\in V_h,
%\end{equation}
%where $C_*$ is (in principle) a computable, $h$ independent, positive constant.
We first consider the homogeneous case, $g\equiv 0$. In this case, we have that
\begin{equation}\label{4.6}
\|u-u_h\|_h\le\|u-E_hu_h\|_h+\|u_h-E_hu_h\|_h,
\end{equation}
and note that the monotonicity of $a_g$ on $H$ implies that
\begin{equation}\label{4.7}
\|u-E_hu_h\|_h^2 = \|D^2(u-E_hu_h)\|_{L^2(\X)}^2\le\frac{a_g(u;u-E_hu_h)-a_g(E_hu_h;u-E_hu_h)}{1-(1-\varepsilon)^\frac{1}{2}}.
\end{equation}
Furthermore, it follows from~(\ref{MTeq}),~(\ref{a:def}), and~(\ref{variation}), that
\begin{equation}\label{4.8}
\begin{aligned}
a_g(u;u-E_hu_h)&-a_g(E_hu_h;u-E_hu_h)\\
&=-(F_\gamma[E_hu_h],\Delta(u-E_hu_h))_{L^2(\X)}\\
&=(-F_\gamma[u_h]+(F_\gamma[u_h]-F_\gamma[E_hu_h]),\Delta(u-E_hu_h))_{L^2(\X)}\\
&\le(\|F_\gamma[u_h]\|_{L^2(\X)}+(\sqrt{d}+1)\|u_h-Eu_h\|_{h})\|u-Eu_h\|_{h},
\end{aligned}
\end{equation}
where we used the inequality 
\begin{equation}\label{16:19}
\sup_{\alpha\in\Lambda}|\gamma^\alpha A^\alpha |\le\sup_{\alpha\in\Lambda}|\gamma^\alpha  A^\alpha -I|+|I|\le\sqrt{1-\varepsilon}+\sqrt{d}<\sqrt{d}+1
\end{equation}
that follows from~(\ref{gamma:2}).
Combining~(\ref{4.7}) and~(\ref{4.8}), we find
\begin{equation}\label{4.9}
\|u-E_hu_h\|_h\le\frac{1}{1-(1-\varepsilon)^\frac{1}{2}}\left(\|F_\gamma[u_h]\|_{L^2(\X)}+(\sqrt{d}+1)\|u_h-E_hu_h\|_h\right),
\end{equation}
which, together with~(\ref{4.6}) implies
\begin{equation}\label{4.10}
\|u-u_h\|_h\le\frac{1}{1-(1-\varepsilon)^\frac{1}{2}}\left(\|F_\gamma[u_h]\|_{L^2(\X)}+(\sqrt{d}+2)\|u_h-E_hu_h\|_h\right).
\end{equation}
In view of~(\ref{14:18}) and~(\ref{4.10}), we arrive at the following \emph{a posteriori} error estimate.
\begin{theorem}\label{thm4.4a}
If $g\equiv 0$, then we have that
\begin{equation}
\|u-u_h\|_h\le\frac{1}{1-(1-\varepsilon)^\frac{1}{2}}\left(\|F_\gamma[u_h]\|_{L^2(\X)}+(\sqrt{d}+2)\sqrt{C_*}\left(\sum_{e\in\Ei}\frac{1}{h_e}\int_e\jump{\partial u_h/\partial n}_e^2\,ds\right)^\frac{1}{2}\right).
\end{equation}
\end{theorem}
We utilise Theorem~\ref{thm4.4a} to prove the analogous result in the inhomogeneous setting.
\begin{theorem}\label{thm4.4}
We have that
\begin{equation}
\begin{aligned}
\|u-u_h\|_h&\le\frac{\|F_\gamma[u_h]\|_{L^2(\X)}+(\sqrt{d}+2)(1+\sqrt{C_*/\sigma})\|g-g_h\|_h}{1-(1-\varepsilon)^\frac{1}{2}}\\
&~~~~~~~~~~~~~~~~~~~~~~~~~~+\frac{(\sqrt{d}+2)\sqrt{C_*}\left(\sum_{e\in\Ei}\frac{1}{h_e}\int_e\jump{\partial u_h/\partial n}_e^2\,ds\right)^\frac{1}{2}}{1-(1-\varepsilon)^\frac{1}{2}}.
\end{aligned}
\end{equation}
\end{theorem}
\emph{Proof:} Define $u_0:=u-g\in H^2(\X)\cap H^1_0(\X)$. We see that $u_0$ satisfies
\begin{align*}
\sup_{\alpha\in\Lambda}\{A^\alpha:D^2u_0- g^\alpha\} =\sup_{\alpha\in\Lambda}\{A^\alpha:D^2u- f^\alpha\} &= 0\quad\mbox{a.e. in }\X,\\
u_0 & = 0\quad\mbox{on }\partial\X,
\end{align*}
where $g^\alpha:=f^\alpha-A^\alpha:D^2g$. Defining $u_{h,0} = u_h-g_h\in V_{h,0}$, by Theorem~\ref{thm4.4a}, we have that
\begin{equation}\label{16:16}
\|u_0-u_{h,0}\|_h\le\frac{\|\sup_{\alpha\in\Lambda}\{\gamma^\alpha(A^\alpha:D^2_hu_{h,0}-g^\alpha)\}\|_{L^2(\X)}+(\sqrt{d}+2)\sqrt{C_*}\left(\sum_{e\in\Ei}\frac{1}{h_e}\int_e\jump{\partial u_{h,0}/\partial n}_e^2\,ds\right)^\frac{1}{2}}{1-(1-\varepsilon)^\frac{1}{2}}.
\end{equation}
Let us denote $g^\alpha_h:=f^\alpha-A^\alpha:D^2_hg_h$. The triangle inequality and~(\ref{16:19}) imply that
\begin{align*}
&\|\sup_{\alpha\in\Lambda}\{\gamma^\alpha(A^\alpha:D^2_hu_{h,0}-g^\alpha)\}\|_{L^2(\X)}\le\|\sup_{\alpha\in\Lambda}\{\gamma^\alpha(A^\alpha:D^2_hu_{h,0}-g_h^\alpha)\}\|_{L^2(\X)}\\
& ~~~~~~~~~~~~~~~+ \|\sup_{\alpha\in\Lambda}\{\gamma^\alpha(A^\alpha:D^2_hu_{h,0}-g^\alpha)\}-\sup_{\alpha\in\Lambda}\{\gamma^\alpha(A^\alpha:D^2_hu_{h,0}-g_h^\alpha)\}\|_{L^2(\X)}\\
&~~~~~~~~~~~~~~~~~~~~~~~~~~~~~~~~~~~~~~~~~~~~~~\le\|F_\gamma[u_h]\|_{L^2(\X)}+(\sqrt{d}+1)\|g-g_h\|_h,
\end{align*}
as well as
$$\left(\sum_{e\in\Ei}\frac{1}{h_e}\int_e\jump{\partial u_{h,0}/\partial n}_e^2\,ds\right)^\frac{1}{2}\le\left(\sum_{e\in\Ei}\frac{1}{h_e}\int_e\jump{\partial u_{h}/\partial n}_e^2\,ds\right)^\frac{1}{2}+\|g-g_h\|_h/\sqrt{\sigma}.$$
Applying the above two estimates to~(\ref{16:16}), and using the triangle inequality once more provides
\begin{align*}
&\|u-u_h\|_h\le\|u_0-u_{h,0}\|_h+\|g-g_h\|_h\\
&~~\le\frac{\|F_\gamma[u_h]\|_{L^2(\X)}+(\sqrt{d}+2)(1+\sqrt{C_*/\sigma})\|g-g_h\|_h+(\sqrt{d}+2)\sqrt{C_*}\left(\sum_{e\in\Ei}\frac{1}{h_e}\int_e\jump{\partial u_h/\partial n}_e^2\,ds\right)^\frac{1}{2}}{1-(1-\varepsilon)^\frac{1}{2}},
\end{align*}
as desired.\hfill$\square$

According to Theorem~\ref{thm4.4}, the error estimator
\begin{equation}\label{4.12}
\begin{aligned}
\eta_h&:=\|F_\gamma[u_h]\|_{L^2(\X)}+\|D^2_h(g-g_h)\|_{L^2(\X)}\\
&~~~~~+\left(\sum_{e\in\Ei}\frac{1}{h_e}\int_e\jump{\partial g_h/\partial n}_e^2\,ds\right)^\frac{1}{2}+\left(\sum_{e\in\Ei}\frac{1}{h_e}\int_e\jump{\partial u_h/\partial n}_e^2\,ds\right)^\frac{1}{2},
\end{aligned}
\end{equation}
is reliable. On the other hand, the local efficiency of $\eta_h$ (modulo data approximation terms) is obvious because
\begin{align}
\|F_\gamma[u_h]\|_{L^2(\X)} \le\|\sup_{\alpha\in\Lambda}\{|\gamma^\alpha A^\alpha:D^2_h(u_h-u)|\}\|_{L^2(\X)}&\le(\sqrt{d}+1)\|D^2_h(u_h-u)\|_{L^2(\X)},\label{4.13}\\
\frac{1}{h_e}\int_e\jump{\partial u_h/\partial n}_e^2\,ds&=\frac{1}{h_e}\int_e\jump{\partial (u-u_h)/\partial n}_e^2\,ds\quad\forall e\in\Ei.\label{4.14}
\end{align}
We denote the local indicators as follows for $e\in\Ei$, and $K\in\Th$:
\begin{align*}
\eta_K(u_h)&:=\|F_\gamma[u_h]\|_{L^2(K)},\quad&\eta_K^{g_h}:=\|D^2(g-g_h)\|_{L^2(K)},\\
\eta_e^2(u_h)&:=\frac{1}{h_e}\|\jump{\partial u_h/\partial n}\|_{L^2(e)}^2,\quad&\eta_e^2(g_h):=\frac{1}{h_e}\|\jump{\partial g_h/\partial n}\|_{L^2(e)}^2.
\end{align*}
\end{section}
\begin{section}{Iterative Scheme}\label{iterative:sec}
Since the form ${a}_h$ is nonlinear in the first argument, we shall employ an iterative scheme, in order to approximate the solution of~(\ref{1})-(\ref{1a}). The method itself is referred to as a semismooth Newton's method (described in Algorithm~\ref{semismooth}), we cannot apply classical Newton's method, since $a_h$ is not classically differentiable in the first argument, due to the presence of the supremum. \begingroup\color{black}The semismooth Newton's method presented is also provided in~\cite{neilan:MT}, and superlinear convergence results for a similar (discontinuous Galerkin) finite element method are proven in~\cite{MR3196952}. This particular semismooth Newton's method is also known as Howard's Algorithm~\cite{bokanowski2009some,howard1960dynamic}.\endgroup

In order to apply the semismooth Newton's method, we iteratively solve discrete problems that correspond to problems of the form~(\ref{2a})--(\ref{2aa}). To this end, given a measurable function $\alpha:\X\to\Lambda$, let us define $a_{\alpha}:V_h\times V_{h,0}\to\mathbb R$, $\ell_\alpha:V_{h,0}\to\mathbb R$ by
\begin{align*}
a_{\alpha}(u,v) &:=(\gamma^\alpha A^\alpha:D^2_hu,\Delta_hv)_{L^2(\X)}+\sum_{e\in\Ei}\frac{\sigma}{h_e}\int_e\jump{\partial u_h/\partial n}_e\jump{\partial v/\partial n}_e\,ds\\
\ell_\alpha(v) &:=(\gamma^\alpha f^\alpha,\Delta_hv)_{L^2(\X)}.
\end{align*}
One can see that the discrete problems of finding $u\in V_h$ such that $u|_{\partial\X} = g_h$, and
$$a_{\alpha}(u,v) = \ell_\alpha(v) \quad\forall v\in V_{h,0},$$
 is equivalent to~(\ref{thescheme2}), in the case that $\Lambda$ is a singleton set.
\begin{algorithm}
\caption{Semismooth Newton's method}\label{semismooth}
 \begin{algorithmic}[1]
\Require{
      $\X\subset\mathbb R^d$,
      $\operatorname{tol}\in\mathbb R^+$,
      $\operatorname{itermax}\in\naturals$,
      $\Th$ a mesh on $\overline{\X}$,
      $V_h$,
      $V_{h,0}$,
      $\Lambda$,
      $\{A^\alpha,\gamma^\alpha,f^\alpha\}_{\alpha\in\Lambda},$
      $u_h^0,g_h\in V_h$
  }
  \State{$k\gets 0$}
  \State{$r\gets 1$}
  \State{$u_h^0\gets u_h^0$}
  \While{$k<\operatorname{itermax}\And r>\operatorname{tol}$}
  \State{Select an arbitrary $\alpha_k\in\operatorname{argmax}F_\gamma[u_h^k]$}
  \State{$u_h^{k+1}\gets\mbox{the solution of}$
  \begin{equation}\label{newdir:lin}
  \begin{aligned}
 a_{\alpha_k}(u,v)  &=\ell_{\alpha_k}(v)\quad\forall v\in V_{h,0},\\
 u|_{\partial\X} & = g_h
 \end{aligned}
  \end{equation}}
  \State{$r\gets\|u_h^{k+1}-u_h^k\|_{L^\infty(\X)}$}
  \State{$u_h^k\gets u_h^{k+1}$}
  \State{$k\gets k+1$}
  \EndWhile
%  \State{return $A_h^{\mathcal{D},k}$}
%   \State{While $k<\operatorname{itermax}$ and $r>\operatorname{tol}$:}
%  \State{Select $\alpha_k\in\Lambda[u^k_h]$}
%  \State{Solve\begin{equation}
%  A_h^{\mathcal{D},k }(u^{k+1}_h,v_h) = \sum_{K\in\Th}\langle\gamma^{\alpha_k}f^{\alpha_k},\Delta v_h\rangle_K\quad\forall v_h\in\dg
%  \end{equation}
%   }
%  \State{$\operatorname{error} = \|u^{k+1}_h-u^k_h\|_{h,1}$}
%  \State{$u^k_h = u^{k+1}_h$}
%  \State{$k=k+1$}
  \end{algorithmic}
  \end{algorithm}
\end{section}

The following algorithm (Algorithm~\ref{AFEM}) describes the adaptive scheme. A general adaptive scheme is defined by iterating the following procedure:
$$\mbox{Solve}\mapsto\mbox{Estimate}\mapsto\mbox{Mark}\mapsto\mbox{Refine}.$$
There are several potential marking schemes that one could consider (for example D\"{o}rfler marking~\cite{dorfler1996convergent}); for the experiments of this section, we implement the maximum marking strategy (described in Algorithm~\ref{AFEM}) with newest vertex bisection (that is, a marked simplex is bisected, and then the generated node is joined to the closest vertex, so that the refinement procedure does not result in hanging nodes).
\begin{algorithm}
\caption{Adaptive finite element method}\label{AFEM}
 \begin{algorithmic}[1]
\Require{
      $\X\subset\mathbb R^d$,
      $\operatorname{tol}\in(0,1)$,
      $\operatorname{itermax}\in\naturals$,
      $\theta\in(0,1]$,
      $\mathscr{T}_{0}$ an initial mesh on $\overline{\X}$
  }
  \State{$k\gets 0,\eta_0\gets1$}
  \While{$k<\operatorname{itermax}\And \eta_k>\operatorname{tol}$}
  \State{\emph{Solve:} $u_{k}\gets\mbox{the solution of Algorithm}~\ref{semismooth}$}
  \State{\emph{Estimate:} For all $K\in\mathscr{T}_k$, $e\in\Ei(\mathscr{T}_k)$, calculate $\eta_K(u_k)$, $\eta_K^{g_k}$, $\eta_e(u_k)$, $\eta_e(g_k)$}
  \State{ $\eta_k\gets\max\{\max_{K\in\mathscr{T}_{k}}\eta_K(u_k),\max_{K\in\mathscr{T}_{k}}\eta_K^{g_k},\max_{e\in\Ei(\mathscr{T}_{k})}\eta_e(u_k),\max_{e\in\Ei(\mathscr{T}_{k})}\eta_e(g_k)\}$}
      \State{\emph{Mark: }} 
\For{$e\in\Ei(\mathscr{T}_k)$}
\If{$\eta_e>\theta\eta_k$}
\State{Mark $e$}
\EndIf
\EndFor
\For{$K\in\mathscr{T}_k$}
\If{$\eta_K>\theta\eta_k$}
\State{Mark $K$}
\EndIf
\EndFor
	  \State{\emph{Refine}: Define $\mathscr{T}_{k+1}$ by bisecting all marked simplices, all simplices whose boundary contains a marked edge, and joining created hanging nodes to closest vertices.}
\State{$k\gets k+1$}
  \EndWhile
%  \State{return $A_h^{\mathcal{D},k}$}
%   \State{While $k<\operatorname{itermax}$ and $r>\operatorname{tol}$:}
%  \State{Select $\alpha_k\in\Lambda[u^k_h]$}
%  \State{Solve\begin{equation}
%  A_h^{\mathcal{D},k }(u^{k+1}_h,v_h) = \sum_{K\in\Th}\langle\gamma^{\alpha_k}f^{\alpha_k},\Delta v_h\rangle_K\quad\forall v_h\in\dg
%  \end{equation}
%   }
%  \State{$\operatorname{error} = \|u^{k+1}_h-u^k_h\|_{h}$}
%  \State{$u^k_h = u^{k+1}_h$}
%  \State{$k=k+1$}
  \end{algorithmic}
  \end{algorithm}
  \subsection{Solving the Linear Problem in FEniCS}
\begingroup\color{black}  At each step of Algorithm~\ref{semismooth}, we are required to solve a linear problem of the form~(\ref{newdir:lin}). This is equivalent to solving a linear system. The following code snippet details how we define the bilinear form $a(\cdot,\cdot)$ and linear form $\ell(\cdot)$ in~(\ref{newdir:lin}) in FEniCS (for simplicity we drop the $\alpha_k$ subscript). For simplicity of exposition, we assume that $A,f,g,\sigma,h$ and $\Th$ are given.
\begin{lstlisting}
# defining finite element
fes = FiniteElement("CG", mesh.ufl_cell(), degree)
    
# defining finite element space
FES = FunctionSpace(mesh,fes)

# defining trial and test functions
uh = Function(FES)
v = TestFunction(FES)

# defining boundary data as L^2 projection
gd = project(g,FES)

# defining unit normal
n = FacetNormal(mesh)

# defining penalty parameter
gamma = (A00+A11)/(pow(A00,2)+2.0*pow(A01,2)+pow(A11,2))

# defining mesh penalty parameter
sig = sigma*pow(h,-1)

# defining jump stabilisation operator
def J_h(u,v,mesh):
    J1 = sig*(n[0]('+')*(u.dx(0)('+')-u.dx(0)('-'))\\
	+ n[1]('+')*(u.dx(1)('+')-u.dx(1)('-')))\\
	*(n[0]('+')*(v.dx(0)('+')-v.dx(0)('-'))\\
	+n[1]('+')*(v.dx(1)('+')-v.dx(1)('-')))\\
	*dS(mesh,metadata={'quadrature_degree': quad_deg})
    return J1

# defining nondivergence part of the bilinear form
def ah(u,v,mesh):
    a = gamma*(A00*u.dx(0).dx(0)+A11*u.dx(1).dx(1)+A01*u.dx(1).dx(0)\\
	+A01*u.dx(0).dx(1))*(v.dx(0).dx(0)+v.dx(1).dx(1))\\
	*dx(mesh,metadata={'quadrature_degree': quad_deg})
    return a

# defining bilinear form a
    a = ah(u,v,mesh)+J_h(u,v,mesh)

# defining linear form l
    l = gamma*(f)*(v.dx(0).dx(0)+v.dx(1).dx(1))\\
	*dx(mesh,metadata={'quadrature_degree': quad_deg})
\end{lstlisting}
\begin{remark}[Boundary data]
We apply the Dirichlet boundary condition using the DirichletBC function in FEniCS. In the code snippet, and in our numerical examples, we take $g_h$ to be the $L^2$ projection of $g$ onto $V_h$. However, one could take $g_h$ to be the unique element of $V_h$ satisifying~(\ref{ghdef}).
\end{remark}\endgroup
\begin{section}{Applications to the Fully Nonlinear Monge--Amp\`{e}re Equation}\label{MA:apps}
%As proven in~\cite{neilan:MT} (cf. Theorem 5), the discrete scheme~(\ref{3.3}) extends to the setting of the fully nonlinear Hamilton--Jacobi--Bellman equation:
%\begin{align}
%\max_{\alpha\in\Lambda}\{A^\alpha:D^2u-f^\alpha\} & = 0,\quad\mbox{a.e. in }\X,\label{HJB}\\
%u & = 0,\quad\mbox{on }\partial\X,\label{HJBa}
%\end{align}
%where $\Lambda$ is a given compact metric space, and then for each $\alpha\in\Lambda$, the functions $A^\alpha,\,f^\alpha$ are defined via $A^\alpha(x):=A(x,\alpha),\,f^\alpha(x):=f(x,\alpha)$, where $A,f\in C(\overline{\X})$. As in the linear case, we must renormalise the equation by multiplying under the maximum in~(\ref{HJB}) by $$\gamma^\alpha(x):=\gamma(x,\alpha)=(A(x,\alpha):I)/(A(x,\alpha):A(x,\alpha)).$$
%One approximates solutions of~(\ref{HJB})-(\ref{HJBa}) by considering the following discrete scheme: to find $u_h\in V_h$ such that
%\begin{equation}\label{HJBdisc}
%\mathscr{A}_h(u_h;v_h):=\int_\X\max_{\alpha\in\Lambda}\{\gamma^\alpha(A^\alpha:D^2_hu-f^\alpha)\}\,\Delta_h v_h+\sum_{e\in\Ei}\frac{\sigma}{h_e}\int_e\jump{\partial u_h/\partial n}_e\jump{\partial v_h/\partial n}_e = 0\quad\forall v_h\in V_h.
%\end{equation}
%In fact, Theorem 5 of~\cite{neilan:MT} shows that the above scheme is uniquely solvable, provided $\sigma$ is sufficiently large. 
%The nonlinear problem~(\ref{1})-(\ref{1a}) arises in many applications of interest, such as finance, engineering, economics, and stochastic optimal control problems. A particular motivation for considering such a problem is in applications to 
Let us consider the fully nonlinear Monge--Amp\`{e}re (MA) equation:
\begin{align}
\operatorname{det}D^2u&= f,\quad\mbox{in}\quad\X,\label{MA-2}\\
u & = g,\quad\mbox{on}\quad\partial\X,\label{MA-1}\\
\begingroup\color{black}u \endgroup&\begingroup\color{black} \quad\mbox{is convex},\endgroup\label{convex10:40}
\end{align}
where $f$ and $g$ are given functions, and $f$ is uniformly positive.
Thanks to~\cite{MR901759} we may characterise equation~(\ref{MA-2})-(\ref{convex10:40}) as the following HJB problem:
\begin{align}
\max_{W\in X}\{-W:D^2u+2f^{1/2}(\operatorname{det}W)^{1/2}\}  & = 0,\quad\mbox{in}\quad\X,\label{HJBMA-1}\\
u & = g,\quad\mbox{on}\quad\partial\X,\label{HJBMA-2}
\end{align}
where $X:=\{W\in\mathbb R^{2\times 2}:W\ge 0,W=W^T,\,\operatorname{Trace}(W)=1\}$. 

However, the control set, $X$, contains degenerate matrices, which do not satisfy~(\ref{3}). This is remedied by the results (in particular Theorem~\ref{2dthm}) below, which prove that we may consider a restricted control set of matrices that satisfy~(\ref{3}) uniformly. \begingroup\color{black}The material that follows is present in~\cite{kawecki:thesis}, under the assumption of classical differentiability of the solution to the MA problem~(\ref{HJBMA-1})-(\ref{HJBMA-2}). Furthermore, similar results are also present in~\cite{MR901759}.\endgroup
 
%The link between the HJB equation and the MA equation allows us to deduce more about the nature of the set $X$. Note that the next lemma holds for any finite dimension $d$.
%\begin{lemma}\label{leq1/d^d}
%The set $X$ satisfies 
%\begin{equation*}
%\sup_{W\in X}\operatorname{det}W=\frac{1}{d^d},
%\end{equation*}
%and a maximiser is given by $\frac{1}{d}I_d$.
%\end{lemma}
%\emph{Proof:} Consider the MA equation 
%\begin{equation*}
%\operatorname{det}D^2u = 1.
%\end{equation*}
%A solution to this problem is given by $u(x)=\frac{1}{2}|x|^2$; note that $D^2u=I_d$. Thus,
%\begin{equation*}
%\begin{aligned}
%0&=\sup_{W\in X}\{-W:D^2u+d(\operatorname{det}W)^{1/d}\}=\sup_{W\in X}\{-\operatorname{Tr}(W)+d(\operatorname{det}W)^{1/d}\}=d\left(\sup_{W\in X}(\operatorname{det}W)^{1/d}\right)-1;
%\end{aligned}
%\end{equation*}
%rearranging, and noting that $a^\frac{1}{d}\ge b$ implies that $a\ge b^d$ for $a,b\ge 0$, we obtain
%\begin{equation*}\sup_{W\in X}\operatorname{det} W=\frac{1}{d^d}.
%\end{equation*}
%It is also clear that $\frac{1}{d}I_d\in X$  and $\operatorname{det}(\frac{1}{d}I_d)=\frac{1}{d^d}$.$\quad\quad\square$
\begin{theorem}\label{wellposedness}
Let $\X$ be a bounded convex open subset of $\mathbb R^2$, and assume that $g\in H^2(\X)$, and that $f\in C^0(\overline{\X})$ is nonnegative. Let $X_\xi:=\{W\in X:\operatorname{det} W\ge\xi\}$. Then, for any constant $\xi\in(0,1/4]$, there exists a unique solution $u\in H^2(\X)$ of the following HJB equation
\begin{equation}
\begin{aligned}
\sup_{W\in X_\xi}\{-W:D^2u+2(\operatorname{det} W)^{1/2}f^{1/2}\}(x)&=0,\quad\mbox{a.e. in }\X,\\
u(x) & = g,\quad\mbox{on }\partial\X.
\end{aligned}
\end{equation}
\end{theorem}
\emph{Proof:} \begingroup\color{black}First note that as $\xi\le1/4$, one has that $(1/2)I_d\in X_\xi$, and so $X_\xi\ne\emptyset$. The set $X_\xi$ also contains only positive definite matrices (since all elements of $X_\xi$ are $2\times 2$ matrices with positive trace and determinant){}{}, and in two dimensions uniform ellipticity implies the Cordes condition. Then, setting $\Lambda=X_\xi$, we can see that $X_\xi$ is a compact metric space; using the Euclidean distance as a metric, and noting that $X_\xi = \mathscr{D}^{-1}([\xi,1/4])$, where $\mathscr{D}:\Lambda\to\mathbb R$ given by
$$\mathscr{D}(W) := \operatorname{det}(W),\quad W\in X_\xi,$$
is a continuous function, we deduce that $X_\xi$ is closed. %Combining the determinant upper bound given by Lemma~\ref{leq1/d^d}, and the unit trace constraint that $X_\xi$ satisfies, we see that $X_\xi$ is also bounded. Thus $X_\xi$ is compact. 
Since each member of $X_\xi$ is of unit trace, denoting the eigenvalues of $W\in X_\xi$ by $\lambda_1,\lambda_2$, we have that $|W|^2=\lambda_1^2+\lambda_2^2=(\lambda_1+\lambda_2)^2-2\lambda_1\lambda_2=1-2\operatorname{det}W\le 1-2\xi<\infty$. Thus $X_\xi$ is bounded. It then follows that $X_\xi$ is compact.\endgroup

We can apply Theorem~\ref{inhomog:eandu:15:11}, yielding existence of a unique $v\in H^2(\X)$ satisfying
\begin{equation}\label{changeofsign:example}
\left\{
\begin{aligned}
\sup_{W\in X_\xi}\{W:D^2v+2(\operatorname{det} W)^{1/2}f^{1/2}\}&=0\,\,\mbox{in}\,\,\X,\\
u&=-g\,\,\mbox{on}\,\,\partial\X.
\end{aligned}
\right.
\end{equation}
We then (uniquely) define $u:=-v$.$\quad\quad\square$
\begin{theorem}\label{2dthm}
Let $d=2$, assume that $\X$ is convex, that $g\in W^{2,\infty}(\X)$, and $f\in C^0(\overline{\X})$ is uniformly positive. Furthermore, assume that $u\in W^{2,\infty}(\X)$ is uniformly convex, and satisfies~(\ref{MA-2})-(\ref{MA-1}). Then, there exists $\xi\in(0,1/4]$ dependent upon $|u|_{W^{2,\infty}(\X)}$, such that $u$ is also the unique solution to 
\begin{equation}\label{17:46}
\left\{
\begin{aligned}
\sup_{W\in X_\xi}\{-W:D^2u+2(\operatorname{det} W)^{1/2}f^{1/2}\}&=0\quad\mbox{a.e. in }\X,\\
u&=g\quad\mbox{on }\partial\X.
\end{aligned}
\right.
\end{equation}
\end{theorem}
\emph{Proof:} Let us define the map $A_u:\overline{\X}\to\mathbb R^{2\times 2}$ by:
\begin{equation}
A_u(x):=\frac{\operatorname{Cof}(D^2u)}{\Delta u},
\end{equation}
note that this map is well defined, since $u$ is uniformly convex, and so, its Laplacian is uniformly positive. Also, since $u\in W^{2,\infty}({\X})$, we have that $A_u\in L^\infty({\X})$. Furthermore, $\operatorname{Cof}(D^2u)$ is symmetric, and 
$$\operatorname{Tr}(A_u) = \frac{1}{\Delta u}\operatorname{Tr}(\operatorname{Cof}(D^2u)) = \frac{\Delta u}{\Delta u}=1,$$
and so $A_u:\overline{\X}\to X$.
%We will now show that $A_u(x)$ is a maximiser of
%\begin{equation}\label{A(x)map:def}
%F_x(W):=-W:D^2u(x)+(\operatorname{det}W)^{1/2}\psi(x),
%\end{equation}
%over $X$ for a.e. $x\in\X$.
%
We see that $A_u$ satisfies
\begin{equation}\label{augiveszero}
\begin{aligned}
&-A_u(x):D^2u(x)+2\operatorname{det}(A_u(x))^{1/2}f^{1/2}\\
&~~ = \frac{1}{\Delta u(x)}(-\operatorname{Cof}(D^2u(x)):D^2u(x)+2(\operatorname{det}(\operatorname{Cof} D^2u(x)))^{1/2}f(x)^{1/2})\\
&~~=\frac{2}{\Delta u(x)}(-\operatorname{det} D^2u(x)+\operatorname{det}(D^2u(x))^{1/2}f(x)^{1/2})\\
&~~=\frac{2}{\Delta u(x)}(-\operatorname{det} D^2u(x)+f(x))=0.
\end{aligned}
\end{equation}
%Note that since $u$ satisfies the MAD problem~(\ref{MA:HJB:sec}), we have that
%$$\max_{W\in X}F_x(W) = 0 = F_x(A_u(x)).$$
%Thus $A_u(x)$ maximises $F_x$ over $X$ for each $x\in\X$.
We also obtain a lower bound on the determinant of $A_u$:
\begin{equation*}
\begin{aligned}
\det(A_u) &= \det\left(\frac{\operatorname{Cof}(D^2u)}{\Delta u}\right)\\
& = \frac{\det(D^2u)}{(\Delta u)^2}\\
& =\frac{f}{(\Delta u)^2}\ge\frac{\delta}{2|u|_{W^{2,\infty}(\X)}^2}=:\xi,
\end{aligned}
\end{equation*}
where $\delta=\inf_{x\in\overline{\X}}f(x)>0$, and so, $\xi>0$. 

Let us consider the following HJB equation: find $v\in H^2(\X)$ such that
\begin{equation}\label{HJBsmaller}
\left\{
\begin{aligned}
\sup_{W\in X_\xi}\{-W:D^2v+2(\operatorname{det} W)^{1/2}f^{1/2}\}&=0,\,\,x\in\X,\\
v&=g,\,\,x\in\partial\X.
\end{aligned}
\right.
\end{equation}
%We also note that since $A_u(x)\in X$ for a.e. $x\in\X$, by Lemma~\ref{leq1/d^d} we have that
%$$\xi\le\operatorname{det}(A(x))\le\frac{1}{4}.$$
There is an important difference between the set $X$ and the set $X_\xi:=\{W\in X:\operatorname{det} W\ge\xi\}$, which is that the latter set consists entirely of positive definite matrices. It then follows from Theorem~\ref{wellposedness} that there exists a unique $v\in H^2(\X)$ that satisfies~(\ref{HJBsmaller}).

We then see that the solution $u$ of the MA equation satisfies (noting that $X_\xi\subseteq X$)
\begin{equation*}
\begin{aligned}
\sup_{W\in X_\xi}\{-W:D^2u+2(\operatorname{det} W)^{1/2}f^{1/2}\}&\le\sup_{W\in X}\{-W:D^2u+2(\operatorname{det} W)^{1/2}f^{1/2}\}=0\quad\mbox{a.e. in }\X.
\end{aligned}
\end{equation*}
Since $A(x)\in X_\xi$ for a.e. $x\in\X$, from~(\ref{augiveszero}), we obtain
\begin{equation*}
\sup_{W\in X_\xi}\{-W:D^2u+2(\operatorname{det} W)^{1/2}f^{1/2}\}\ge-A(x):D^2u+2(\operatorname{det} A(x))^{1/2}f^{1/2}=0\quad\mbox{a.e. in }\X.
\end{equation*}
By combining these results, we obtain
\begin{equation*}
\sup_{W\in X_\xi}\{-W:D^2u(x)+2(\operatorname{det} W)^{1/2}
f^{1/2}\}=0\quad\mbox{a.e. in }\X.
\end{equation*}
Since $u=g$ on $\partial\X$, and $u\in H^2(\X)$, by uniqueness $u=v$.$\quad\quad\square$
% Furthermore, if $u_1,u_2\in W^{2,\infty}(\X)$, are two uniformly convex functions satisfying~(\ref{17:46}), then $u_1,u_2\in W^{2,\infty}(\X)\subset H^2(\X)$, and they both satisfy~(\ref{HJBsmaller}), with $\xi=\xi^*$, where
%$$\xi^*:=\min\left\{\frac{\delta}{2|u_1|_{W^{2,\infty}(\X)}},\frac{\delta}{2|u_2|_{W^{2,\infty}(\X)}}\right\}.$$ Thus, by uniqueness, $u_1\equiv u_2$, that is,~(\ref{17:46}) is also uniquely solvable in the class of uniformly convex $W^{2,\infty}(\X)$ functions.$\quad\quad\square$
\end{section}
\begin{section}{Numerical Results}\label{sec:5}
   \begin{remark}[PDE coefficients]
%In Experiments~\ref{exp1} and~\ref{exp3}, we consider the coefficient matrix given by $A_{ij}:=(1+\delta_{ij})\frac{x_ix_j}{|x_i||x_j|}$ composed with an affine map in Experiment~\ref{exp1}, and a nonaffine map in Experiment~\ref{exp3}. This example was considered in~\cite{MR3077903}, and the composition with a nonaffine map was considered in~\cite{gallistl2017variational}. Furthermore, in Experiment~\ref{exp1} we multiply the coefficient matrix by an interface function $\chi_\X$ (defined below), so that the coefficients have large jumps (this was inspired by the elliptic interface experiment in~\cite{}). 
In Experiment~\ref{exp1}, we consider the coefficient matrix given by $A_{ij}:=(1+\delta_{ij})\frac{x_ix_j}{|x_i||x_j|}$ composed with an affine map. This example was considered in~\cite{MR3077903}. Furthermore, we multiply the coefficient matrix by an interface function $\chi_\X$ (defined below), so that the coefficients have large jumps.
\end{remark}
% (this was inspired by the elliptic interface experiment in~\cite{}). 
\begin{remark}[Monge--Amp\`{e}re]
In Experiment~\ref{exp4}, we consider a family of Monge--Amp\`{e}re type problems with true solutions that have been slightly modified from an example that is present in~\cite{MR3162358} (cf.~\cite{MR3162358}, Test 4). The modifications ensure that the true solutions are uniformly convex and belong to $W^{2,\infty}(\X)\setminus V_h$.
\end{remark}
\subsection{Experiment 1}\label{exp1}
In this experiment, we consider the following problems
\begin{equation}\label{exp:1-2}
\left\{
\begin{aligned}
\sum_{i,j=1}^2(1+\delta_{ij})\frac{(x_i-0.5)}{|x_i-0.5|}\frac{(x_j-0.5)}{|x_j-0.5|} \chi^N_\X(x_1,x_2)D^2_{ij}u_s&= f_s,\quad\mbox{in}\quad\X,\\
u_s & = g_s,\quad\mbox{on}\quad\partial\X,
\end{aligned}
\right.
\end{equation}where $\X=(0,1)^2$. Furthermore, the interface function $\chi^N_\X$ satisfies $\chi^N_\X = 1$ on $\X_1:=\cup_{i,j=0}^{N/2-1}\{2i/N<x_1<(2i+1)/N,2j/N<x_2<(2j+1)/N\}$, and $\chi^N_\X = 1000$ on $\X\setminus\X_1$. In this case we take $N=20$.
%$\{x_1,x_2<0.5\}\cup\{x_1,x_2>0.5\}$, $\chi_\X = 100$ on $\{x_1<0.5,x_2>0.5\}$, and $\chi_\X = 200$ on $\{x_1>0.5,y_1<0.5\}$.
In this case $f_s$ and $g_s$ are chosen so that the solution of~(\ref{exp:1-2}) is given by $u(x) = |x|^{1+s}.$ We consider the exponent $s\in\{0.01,0.1,0.2,\ldots,0.5\}$. It holds that $u_s\in H^{2+\delta}(\X)$, for arbitrary $\delta\in[0,s]$. Furthermore, $u_s$ lacks regularity at the origin, and one can see in Figure~\ref{plot:exp1}, the error estimator prioritises refinement towards the origin, in addition to further refinement in the areas of the domain where $\chi_\X$ is the largest.
We apply both Algorithm~\ref{AFEM} with $\theta = 0.2$, and a uniform refinement procedure, so that we may compare the two approaches. For clarity, we denote the numerical solution by $u_{h,\operatorname{adapt}}$, and $u_{h,\operatorname{unif}}$ for the adaptive and uniform approach, respectively. We consider a variety of values of $s$, and polynomial degree, $p$, and calculate the error in the following (semi) norms: $\|\cdot\|_{L^2(\X)}$, $|\cdot|_{H^1(\X)}$, $\|\cdot\|_{h}$, and also calculate the error estimator $\eta_h$.

\textbf{Case 1:} $p=4$, and $s=0.01$. We observe that 
\begin{align*}\|u_{0.01}-u_{h,\operatorname{adapt}}\|_{L^2(\X)}&=\mathcal{O}(\operatorname{ndofs}^{-2}),\quad\|u_{0.01}-u_{h,\operatorname{unif}}\|_{L^2(\X)}=\mathcal{O}(\operatorname{ndofs}^{-1.01}),\\
|u_{0.01}-u_{h,\operatorname{adapt}}|_{H^1(\X)}&=\mathcal{O}(\operatorname{ndofs}^{-1}),\quad|u_{0.01}-u_{h,\operatorname{unif}}|_{H^1(\X)}=\mathcal{O}(\operatorname{ndofs}^{-0.51}),\\
\eta_{\operatorname{adapt}},\|u_{0.01}-u_{h,\operatorname{adapt}}\|_{h}&=\mathcal{O}(\operatorname{ndofs}^{-0.01}),\quad\eta_{\operatorname{unif}},\|u_{0.01}-u_{h,\operatorname{unif}}\|_h=\mathcal{O}(\operatorname{ndofs}^{-0.005}),
\end{align*}
and so, the adaptive method outperforms the uniform scheme. We also plot the effectivity index in Figure~\ref{plot:exp2}, verifying~(\ref{4.13})--(\ref{4.14}), for the adaptive scheme. %The exact values are also provided in Tables~\ref{table1}-\ref{table2}.

\textbf{Case 2:} $p=3$, and $s\in\{0.1,\ldots,0.5\}$.
%We plot the error estimator, $\eta_h$, and the errors $\|u_{0.01}-u_{h,0.01}\|_{h}$, $|u_{0.01}-u_{h,0.01}|_{H^1(\X)}$, and $\|u_{0.01}-u_{h,0.01}\|_{L^2(\X)}$, in Figure~\ref{plot:exp1}, where $u_h$ is the solution to the adaptive scheme $p=4$. We observe convergence of the estimator, $\eta_h$ as well as the true error $\|u_{0.01}-u_{h,0.01}\|_{h}$ at the order $\mathcal{O}(\operatorname{ndofs}^{-0.01})$, due to the lack of regularity of the true solution. In contrast, we see much faster convergence in the $|\cdot|_{H^1(\X)}$-seminorm and $\|\cdot\|_{L^2(\X)}$-norm. In particular, we observe that $|u_{0.01}-u_{h,0.01}|_{H^1(\X)}=\mathcal{O}(\operatorname{ndofs}^{-1})$ and $\|u_{0.01}-u_{h,0.01}\|_{L^2(\X)}=\mathcal{O}(\operatorname{ndofs}^{-2})$.
We observe that 
\begin{align*}\|u_{s}-u_{h,\operatorname{adapt},s}\|_{L^2(\X)}&=\mathcal{O}(\operatorname{ndofs}^{-(2+s)}),\quad&\|u_{0.01}-u_{h,\operatorname{unif},s}\|_{L^2(\X)}=\mathcal{O}(\operatorname{ndofs}^{-(1+s/2)}),\\
|u_{s}-u_{h,\operatorname{adapt},s}|_{H^1(\X)}&=\mathcal{O}(\operatorname{ndofs}^{-(1+s)}),\quad&|u_{0.01}-u_{h,\operatorname{unif},s}|_{H^1(\X)}=\mathcal{O}(\operatorname{ndofs}^{-(0.5+s/2)}),\\
\|u_{s}-u_{h,\operatorname{adapt},s}\|_{h}&=\mathcal{O}(\operatorname{ndofs}^{-s}),\quad&\|u_{0.01}-u_{h,\operatorname{unif},s}\|_h=\mathcal{O}(\operatorname{ndofs}^{-s/2}),\\
\eta_{\operatorname{adapt},s}&=\mathcal{O}(\operatorname{ndofs}^{-1.11}),\quad&\eta_{\operatorname{unif},s}=\mathcal{O}(\operatorname{ndofs}^{-s/2}).
\end{align*}
%
%We observe an improvement in convergence rates for both the estimator, $\eta_h$, and for the error $\|u_s-u_{h,s,\operatorname{adapt}}\|_{h}$ as the parameter $s$ increases, these are all plotted in Figure~\ref{plot:exp2}. In particular, we observe the convergence rates $\|u_s-u_{h,s,\operatorname{adapt}}\|_{h}=\mathcal{O}(\operatorname{ndofs}^{-s})$. 
%Furthermore, the exact convergence rates in $\|\cdot\|_{L^2(\X)}$, $|\cdot|_{H^1(\X)}$ (as well as $\eta_h$, and $\|\cdot\|_h$) are provided in Tables~\ref{table3}-\ref{table12}. 
%Again, we see that the adaptive scheme outperforms the uniform refinement scheme.
\begin{figure}[h]
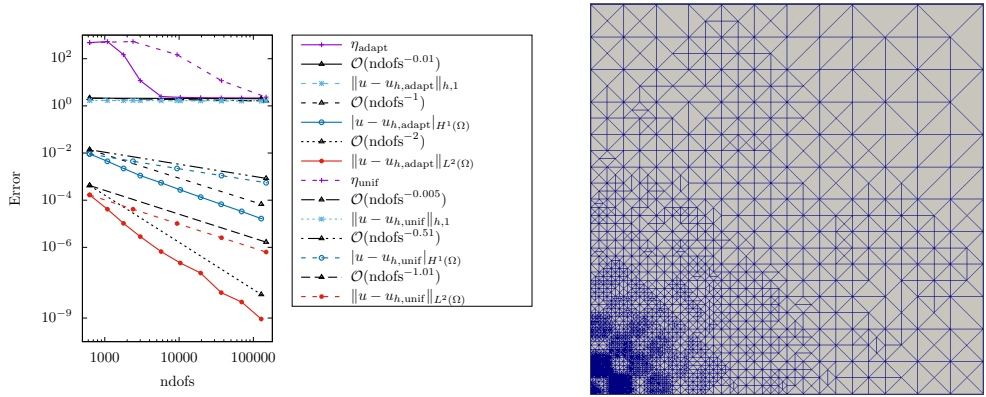

\begin{center}
  \begin{tabular}{lrr}
  \hspace{-3cm}
    \scalebox{0.56}{%
    \renewenvironment{table}[1][]{\ignorespaces}{\unskip}%
    \input{new_paper_exp_0_01.tex}
    \unskip
    }% P1 norec
    &
    \hspace{0.1cm}
   %\vspace{20cm}
    \scalebox{0.23}{%
    \renewenvironment{table}[1][]{\ignorespaces}{\unskip}%
   \input{exp_1_refinement.tex}
    \unskip
    }% P1 norec
  \end{tabular}
  \end{center}
 \caption{On the left are the convergence rates for Experiment~\ref{exp1}, with $s=0.01$, and on the right is the final adapted mesh.}
\label{plot:exp1}
\end{figure}
\begin{figure}[h]
\begin{center}
%  \begin{tabular}{lrr}
%  \hspace{-1.7cm}
%    \scalebox{0.7}{%
%    \renewenvironment{table}[1][]{\ignorespaces}{\unskip}%
%    \input{new_paper_exp_1_reg_tests.tex}
%    \unskip
%    }% P1 norec
%    &
%    \hspace{-1cm}
%   %\vspace{20cm}
    \scalebox{0.7}{%
   \input{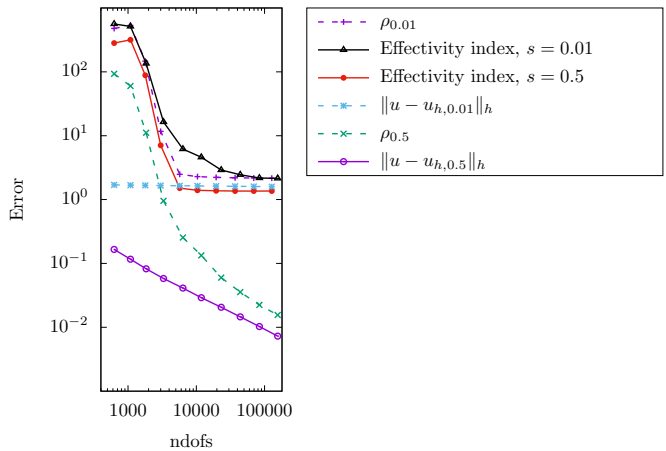}}
%   \unskip
%    }% P1 norec
%  \end{tabular}
  \end{center}
 \caption{%Convergence rates for Experiment~\ref{exp1}, for $s\in\{0.1,\ldots,0.5\}$, with polynomial degree $p=3$ (left). 
 Plot of the effectivity index, with error indicators and true errors for Experiment~\ref{exp1}, with polynomial degree $p=4$.}
\label{plot:exp2}
\end{figure}
\newpage
\subsection{Experiment 2}\label{exp2}
In the previous experiment, we observed the advantage of applying the adaptive scheme, when compared with uniform refinement (see Figure~\ref{plot:exp1}). However, the solution is known to possess $H^s$-regularity, with $s>2$, and is known to lack regularity at the origin. We propose a second experiment, in which the solution is unknown, the right-hand is smooth, and the coefficients are indeed discontinuous (we choose a smooth right-hand side in order to surmise that any bad behaviour of the solution is due to the coefficients and regularity of $\partial\X$). In particular we consider the boundary value problem:
\begin{equation}\label{exp:2:2}
\left\{
\begin{aligned}
\sum_{i,j=1}^2(1+\delta_{ij})\frac{(x_i-0.5)}{|x_i-0.5|}\frac{(x_j-0.5)}{|x_j-0.5|} \chi^N_\X(x_1,x_2)D^2_{ij}u_N&= 1,\quad\mbox{in}\quad\X,\\
u_N & = 0,\quad\mbox{on}\quad\partial\X,
\end{aligned}
\right.
\end{equation}where $\X=(0,1)^2$. We consider the case $N=10$, and in this case the PDE theory implies that $u_N\in H^2(\X)\cap H^1_0(\X)$ (see~(\ref{variation})). We consider the polynomial degree $p=2$, and an initial triangulation with a resolution that matches the indicator function (i.e., $N$ squares in each coordinate direction, with each square further bisected into two triangles), and apply uniform mesh refinement, as well as adapted refinement (applying Algorithm~\ref{AFEM}), and compare the results. 

The solution is unknown, and so we plot the error estimator $\eta_h$ in each case. Due to discrete Poincar\'{e}--Friedrichs' inequalities and the reliability and efficiency of the estimator, $\eta_h$ may be used to as a predictor for the (semi)norms $\|\cdot\|_{L^2(\X)}$, $|\cdot|_{H^1(\X)}$, and $\|\cdot\|_h$. Since the convergence rates in $\|\cdot\|_{L^2(\X)}$, $|\cdot|_{H^1(\X)}$, as predicted by $\eta_h$ are likely to be pessimistic, we calculate the error arising between successive meshes, and appeal to this to guide the convergence. In particular, we define $\theta_k:=u_k-u_{k-1}$, where the subscript $k$ denotes the current refinement level, and appeal to the fact that for the norms under consideration $\|u-u_{k}\|\le\|u-u_{k-1}\|+\|\theta_k\|$, and that the contribution $\|\theta_k\|$ should be the dominating term. We plot the final adapted mesh generated by the adaptive scheme in Figure~\ref{plot:exp:2:2}. The predictions show that the adaptive scheme outperforms uniform refinement, however, not to the same degree as is observed in Experiment~\ref{exp1}, in the $L^2$- and $H^1$-norms. The exact values are provided in Tables~\ref{table13}-\ref{table15}.
%\begin{figure}
%%\begin{center}
%%  \begin{tabular}{lrr}
%  \hspace{-3cm}
% % \vspace{3cm}
%%    \scalebox{0.5}{%
%%    \renewenvironment{table}[1][]{\ignorespaces}{\unskip}%
%%    \input{new_paper_experiment_2.tex}
%%    \unskip
%%    }% P1 norec
%%    &
%   % \hspace{0.1cm}
%   %\vspace{20cm}
%%
%%    \renewenvironment{table}[1][]{\ignorespaces}{\unskip}%
%    \scalebox{0.2}{\input{exp_2_new_refinement_N10.tex}}
%%    \unskip
%%    }% P1 norec
%%  \end{tabular}
%%  \end{center}
%% \caption{On the left are the convergence rates for Experiment~\ref{exp:2:2}, with $N=10$, and on the right is the final adapted mesh.}
%
%\label{plot:exp:2:2}
%\end{figure}
\vspace{2.5cm}
\begin{figure}
\begin{center}
\hspace{-2.5cm}
\scalebox{0.2}{\input{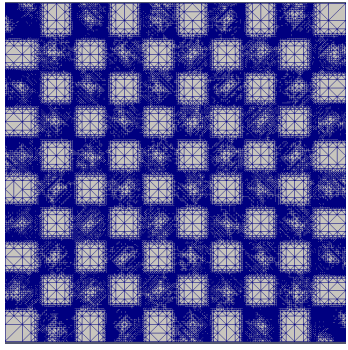}}
\end{center}
\caption{Final adapted mesh for Experiment~\ref{exp:2:2}, with $N=10$.}\label{plot:exp:2:2}
\end{figure}
%%% uniform
\begin{table}[h]
%\centering
 \caption{$p=2$, uniform refinement}\label{table13}
\begin{subtable}{1\linewidth}
 \resizebox{\linewidth}{!}{\begin{filecontents*}{exp2_table-unif-ndofinc_unknown.csv}
1.681000000000000000e+03 2.122458750246703496e-01 0.000000000000000000e+00 3.548413361775816597e-03 0.000000000000000000e+00 4.531484763742778124e-04 0.000000000000000000e+00 1.695052919058331042e+02 0.000000000000000000e+00
6.561000000000000000e+03 2.519516632709564163e-01 1.259345256558220827e-01 3.845563896771360055e-03 5.905594970694931500e-02 1.549213452701842994e-04 -7.881762471593988240e-01 1.194128059019715806e+02 -2.572400407422775204e-01
2.592100000000000000e+04 1.613406043505953569e-01 -3.244167630463026253e-01 1.389417834848707739e-03 -7.409765575096584067e-01 4.709122112266291173e-05 -8.667456951933036891e-01 7.743279863485645365e+01 -3.152869321070880249e-01
1.030410000000000000e+05 9.638974229271275729e-02 -3.732539542554948464e-01 4.209872378828287236e-04 -8.651985427176460686e-01 1.270048228359197777e-05 -9.495484206438333796e-01 4.545889447355671109e+01 -3.859228802356640120e-01
\end{filecontents*}
\pgfplotstabletypeset[dec sep align,
   fixed zerofill,
   precision=3,
   col sep=space,
%      columns/0/.style ={column name=$h$},
   columns/0/.style ={int detect, column name=$\operatorname{ndofs}$},
   columns/1/.style ={column name=$\|\theta_h\|_h$},
   columns/2/.style ={column name=$\operatorname{EOC}$},
        columns/3/.style ={column name=$|\theta_h|_{H^1(\Omega)}$},
   columns/4/.style ={column name=$\operatorname{EOC}$},
   columns/5/.style ={column name=$\|\theta_h\|_{L^2(\Omega)}$},
   columns/6/.style ={column name=$\operatorname{EOC}$},
         columns/7/.style ={column name=$\eta_h$},
   columns/8/.style ={column name=$\operatorname{EOC}$},
%   columns/10/.style ={column name=},
   every head row/.style={before row=\toprule,after row=\midrule},
every last row/.style={after row=\bottomrule}
   ]{exp2_table-unif-ndofinc_unknown.csv}}
%\hspace{-1cm}
%\centering
\end{subtable}
\bigskip
\caption{$p=2$, adaptive refinement}\label{table15}
\begin{subtable}{1\linewidth}
 \resizebox{\linewidth}{!}{
\begin{filecontents*}{exp2_table-adapt_unknown-inc.csv}
1.093000000000000000e+03 1.864728018995168435e-01 0.000000000000000000e+00 3.125903829273612641e-03 0.000000000000000000e+00 4.041219076573664918e-04 0.000000000000000000e+00 2.342294214695764367e+02 0.000000000000000000e+00
2.911000000000000000e+03 1.880660817413866426e-01 8.685445252804015540e-03 2.923799376628850432e-03 -6.823351956848433542e-02 2.728930655823778710e-04 -4.008252439134215028e-01 1.722884840212084043e+02 -3.135361651836178210e-01
7.501000000000000000e+03 1.889315667301778578e-01 4.850794005372553383e-03 2.290898258690540756e-03 -2.577175987734097173e-01 5.138842271548416097e-05 -1.763969481764116720e+00 1.270457873539433393e+02 -3.218277680759111248e-01
2.015900000000000000e+04 1.366672515961921963e-01 -3.275652134719881370e-01 1.285465197918719349e-03 -5.844778954953223638e-01 5.096939950178126512e-05 -8.281757213219035144e-03 8.203369203180798763e+01 -4.424550856315543901e-01
5.250300000000000000e+04 8.344961956706273276e-02 -5.153531390932584166e-01 4.466355411431319547e-04 -1.104378950279725435e+00 2.723899950557155222e-05 -6.545789230720853169e-01 5.096777995544849205e+01 -4.972071585739025190e-01
1.329730000000000000e+05 5.011229683774671290e-02 -5.487893619255238553e-01 1.470176244262275946e-04 -1.195759597980489586e+00 8.606970541486340590e-06 -1.239758274128927251e+00 3.166640193944383341e+01 -5.121595008725178255e-01
\end{filecontents*}
\pgfplotstabletypeset[dec sep align,
   fixed zerofill,
   precision=3,
   col sep=space,
%      columns/0/.style ={column name=$h$},
   columns/0/.style ={int detect, column name=$\operatorname{ndofs}$},
   columns/1/.style ={column name=$\|\theta_h\|_h$},
   columns/2/.style ={column name=$\operatorname{EOC}$},
        columns/3/.style ={column name=$|\theta_h|_{H^1(\Omega)}$},
   columns/4/.style ={column name=$\operatorname{EOC}$},
   columns/5/.style ={column name=$\|\theta_h\|_{L^2(\Omega)}$},
   columns/6/.style ={column name=$\operatorname{EOC}$},
         columns/7/.style ={column name=$\eta_h$},
   columns/8/.style ={column name=$\operatorname{EOC}$},
%   columns/10/.style ={column name=},
   every head row/.style={before row=\toprule,after row=\midrule},
every last row/.style={after row=\bottomrule}
   ]{exp2_table-adapt_unknown-inc.csv}}
    \end{subtable}
    \end{table}
\subsection{Experiment 3}\label{exp4}
In this experiment, we consider the following Monge--Amp\`{e}re problems
\begin{align}
\operatorname{det}D^2u_a&= f_a,\quad\mbox{in}\quad\X,\label{MA1}\\
u_a & = g_a,\quad\mbox{on}\quad\partial\X\label{MA2},
\end{align}
on $\X=(0,1)^2$. %Thanks to~\cite{MR901759} we may characterise equation~(\ref{MA1}) as the following HJB problem:
%$$\max_{W\in X_\xi}\{-W:D^2u_a+2f_a^{1/2}(\operatorname{det}W)^{1/2}\}  = 0,$$
%where $X_\xi:=\{W\in X:\operatorname{det}W\ge\xi\}$, and $X:=\{W\in\mathbb R^{2\times 2}:W^T=W,\,\operatorname{Trace}(W)=1\}$ (note that we take $\xi=1/100$).

\textbf{Case 1:} The functions $f_a$ and $g_a$ are chosen so that the true solution of~(\ref{MA1})-(\ref{MA2}) is given by
$$u_a(x_1,x_2) = |x_1-a|\sin(x_1-a)+50.0(x_1^2+x_2^2),$$
for $a\in\{0.4,0.5\}$. Our initial mesh is a uniform triangulation on $\overline{\X}$ consisting of two squares (each further subdivided into two right-angled triangles) in the $x_1$ and $x_2$ direction, as such, we have that $u_{0.5}$ is piecewise smooth on the initial mesh (and all subsequent meshes, since each marked triangle is bisected), however, $u_{0.4}$ does not enjoy this piecewise smoothness property, and so, its approximation, $u_{h,a=0.4}$, does not converge as fast, as observed in Figure~\ref{plot:exp4:1}. In both cases, $u_a\in W^{2,\infty}(\X)$, and we set the polynomial degree $p=4$. Note that in this case we apply the adaptive finite element method given by Algorithm~\ref{AFEM}, in conjunction with the semismooth Newton's method given by Algorithm~\ref{semismooth}.

We also compare the adaptive scheme with that of uniform refinement. The exact results are provided in Tables~\ref{table16}-\ref{table21}. 
We observe that when $a=0.4$, the adaptive scheme out performs uniform refinement, whereas when $a=0.5$ the two approaches are comparable (we surmise this is due to the piecewise smoothness property of $u_{0.5}$).

\textbf{Case 2:} Here we take $f_a\equiv1$, $g_a\equiv 0$. In this case the true solution is unknown, and so we rely on the error estimator, as well as the the incremental solutions in order to indicate the performance of the numerical method (as in Experiment~\ref{exp2}). We take $p=4$ and compare the adaptive scheme with uniform refinement. We display the exact convergence results in Tables~\ref{table22}-\ref{table24}, and observe that the adaptive scheme outperforms the uniform scheme in all (semi)norms.
%%% s = 0.01, uniform
%%%% s = 0.1, adapt
%\begin{table}
%\hspace{-1cm}
%%\centering
%\pgfplotstabletypeset[dec sep align,
%   fixed zerofill,
%   precision=3,
%   col sep=space,
%%      columns/0/.style ={column name=$h$},
%   columns/0/.style ={int detect, column name=$h$},
%   columns/1/.style ={column name=$\|e_h\|_h$},
%   columns/2/.style ={column name=$\operatorname{EOC}$},
%        columns/3/.style ={column name=$|e_h|_{H^1(\Omega)}$},
%   columns/4/.style ={column name=$\operatorname{EOC}$},
%   columns/5/.style ={column name=$\|e_h\|_{L^2(\Omega)}$},
%   columns/6/.style ={column name=$\operatorname{EOC}$},
%         columns/7/.style ={column name=$\eta_h$},
%   columns/8/.style ={column name=$\operatorname{EOC}$},
%%   columns/10/.style ={column name=},
%   every head row/.style={before row=\toprule,after row=\midrule},
%every last row/.style={after row=\bottomrule}
%   ]{data/table-unif_MA_05.csv}
%    \caption{$p=4$, $a=0.5$, uniform refinement (mesh size EOCs)}
%\end{table}
%
%%% s = 0.3, uniform
\begin{table}
\caption{$p=4$, $a=0.4$, uniform refinement}\label{table16}
\begin{subtable}{1\linewidth}
%\centering
 \resizebox{\linewidth}{!}{
\begin{filecontents*}{table-unif-ndof_MA_04.csv}
8.100000000000000000e+01 8.319673549676219748e-01 0.000000000000000000e+00 2.220017590911668379e-02 0.000000000000000000e+00 1.476509099912783323e-03 0.000000000000000000e+00 8.795067434795964534e-01 0.000000000000000000e+00
2.890000000000000000e+02 7.027967599567351842e-01 -1.326481425088571753e-01 9.744669719462866961e-03 -6.473225752474800476e-01 3.105287788650058756e-04 -1.225775802406688486e+00 7.973065571018659314e-01 -7.714129012633685889e-02
1.089000000000000000e+03 3.229047032654048244e-01 -5.862485171252111993e-01 2.366406700748896995e-03 -1.066907986181331758e+00 3.585623831720451593e-05 -1.627301250655978304e+00 3.586394902346060065e-01 -6.022376979752466220e-01
4.225000000000000000e+03 2.676734194740749140e-01 -1.383649544651049668e-01 1.211022511727656743e-03 -4.941198139530580957e-01 1.231020048278880637e-05 -7.885538717374918072e-01 3.328655939366133576e-01 -5.500894247378517593e-02
1.664100000000000000e+04 1.423274146126781436e-01 -4.607633142132090409e-01 2.989447159594469912e-04 -1.020506545913359098e+00 9.441451260278993125e-06 -1.935430358145217300e-01 1.635059533449582103e-01 -5.185755156723906856e-01
\end{filecontents*}
\pgfplotstabletypeset[dec sep align,
   fixed zerofill,
   precision=3,
   col sep=space,
%      columns/0/.style ={column name=$h$},
   columns/0/.style ={int detect, column name=$\operatorname{ndofs}$},
   columns/1/.style ={column name=$\|e_h\|_h$},
   columns/2/.style ={column name=$\operatorname{EOC}$},
        columns/3/.style ={column name=$|e_h|_{H^1(\Omega)}$},
   columns/4/.style ={column name=$\operatorname{EOC}$},
   columns/5/.style ={column name=$\|e_h\|_{L^2(\Omega)}$},
   columns/6/.style ={column name=$\operatorname{EOC}$},
         columns/7/.style ={column name=$\eta_h$},
   columns/8/.style ={column name=$\operatorname{EOC}$},
%   columns/10/.style ={column name=},
   every head row/.style={before row=\toprule,after row=\midrule},
every last row/.style={after row=\bottomrule}
   ]{table-unif-ndof_MA_04.csv}}
\end{subtable}
\bigskip
\caption{$p=4$, $a=0.4$, adaptive refinement}
\begin{subtable}{1\linewidth}
%%% s = 0.3, adapt
%\centering
 \resizebox{\linewidth}{!}{
\begin{filecontents*}{table-adapt_MA_04.csv}
8.100000000000000000e+01 8.319673549676219748e-01 0.000000000000000000e+00 2.220017590911668379e-02 0.000000000000000000e+00 1.476509099912783323e-03 0.000000000000000000e+00 8.795067434795964534e-01 0.000000000000000000e+00
1.390000000000000000e+02 8.875259736955205536e-01 1.197066994000138412e-01 1.738296465171706648e-02 -4.529598248050812903e-01 9.446175383099503524e-04 -8.271023070425088664e-01 9.525981681294961589e-01 1.478301310939861835e-01
3.410000000000000000e+02 6.097133580470659453e-01 -4.183700362747865831e-01 9.273286013506045616e-03 -7.001859885390304683e-01 3.332560382789132784e-04 -1.160975179848680927e+00 7.120695884240395568e-01 -3.242865515964268308e-01
8.390000000000000000e+02 4.682266719248665554e-01 -2.932668475186474000e-01 3.158708186189492732e-03 -1.196202249292711572e+00 7.298733970237291833e-05 -1.686746046233475127e+00 5.153088363927385007e-01 -3.592125954306355506e-01
1.961000000000000000e+03 2.600185043115390560e-01 -6.928154668257873183e-01 1.426235348323863983e-03 -9.365437188636835941e-01 2.890805809333661742e-05 -1.090891136207853496e+00 3.173634373374746853e-01 -5.709296803672970766e-01
4.385000000000000000e+03 1.844031679336662621e-01 -4.270079898182347344e-01 5.369324134193491577e-04 -1.213966468208026450e+00 2.487634340782225953e-05 -1.866491210262099520e-01 2.174459221387545771e-01 -4.698408600091931175e-01
9.341000000000000000e+03 1.171783461716435909e-01 -5.995942708067188986e-01 1.815020834060908426e-04 -1.434238457645230724e+00 3.954449002212071086e-06 -2.431920666636696815e+00 1.482641964849404359e-01 -5.064035927303094375e-01
1.960900000000000000e+04 8.074002034444287268e-02 -5.022587957835579209e-01 7.138714128618598353e-05 -1.258333892145410360e+00 3.293469328147044934e-06 -2.466367889308459238e-01 9.936527178720834519e-02 -5.396526859057430014e-01
\end{filecontents*}
\pgfplotstabletypeset[dec sep align,
   fixed zerofill,
   precision=3,
   col sep=space,
%      columns/0/.style ={column name=$h$},
   columns/0/.style ={int detect, column name=$\operatorname{ndofs}$},
   columns/1/.style ={column name=$\|e_h\|_h$},
   columns/2/.style ={column name=$\operatorname{EOC}$},
        columns/3/.style ={column name=$|e_h|_{H^1(\Omega)}$},
   columns/4/.style ={column name=$\operatorname{EOC}$},
   columns/5/.style ={column name=$\|e_h\|_{L^2(\Omega)}$},
   columns/6/.style ={column name=$\operatorname{EOC}$},
         columns/7/.style ={column name=$\eta_h$},
   columns/8/.style ={column name=$\operatorname{EOC}$},
%   columns/10/.style ={column name=},
   every head row/.style={before row=\toprule,after row=\midrule},
every last row/.style={after row=\bottomrule}
   ]{table-adapt_MA_04.csv}}
\end{subtable}
\bigskip
%%% s = 0.2, uniform
\caption{$p=4$, $a=0.5$, uniform refinement}
\begin{subtable}{1\linewidth}
%\centering
 \resizebox{\linewidth}{!}{
\begin{filecontents*}{table-unif-ndof_MA_05.csv}
8.100000000000000000e+01 8.885535614132523009e-04 0.000000000000000000e+00 3.135685197569323756e-05 0.000000000000000000e+00 2.558004085498390167e-06 0.000000000000000000e+00 9.901844463934765094e-04 0.000000000000000000e+00
2.890000000000000000e+02 9.264761541178043723e-05 -1.777383369203534791e+00 2.092001425253930763e-06 -2.128427199656077740e+00 8.780566070035832300e-08 -2.650877487538685795e+00 1.082114470267785475e-04 -1.740442752465964471e+00
1.089000000000000000e+03 1.114214818135774922e-05 -1.596627943037243202e+00 1.308236070117833511e-07 -2.089590508428970317e+00 2.807900142266242850e-09 -2.595144776787950747e+00 1.331528618337638479e-05 -1.579370352215061279e+00
4.225000000000000000e+03 1.319624717782388893e-06 -1.573573919433395085e+00 8.141733147586514797e-09 -2.048185537349811902e+00 8.958509098272433380e-11 -2.541013729243995378e+00 1.596974469581623588e-06 -1.564290687087396714e+00
1.664100000000000000e+04 1.672644227823031959e-07 -1.506748612227992057e+00 5.777111223242010101e-10 -1.929958696669161311e+00 5.610671244444278208e-11 -3.413454150838874379e-01 2.055690140241296431e-07 -1.495483785384092101e+00
\end{filecontents*}
\pgfplotstabletypeset[dec sep align,
   fixed zerofill,
   precision=3,
   col sep=space,
%      columns/0/.style ={column name=$h$},
   columns/0/.style ={int detect, column name=$\operatorname{ndofs}$},
   columns/1/.style ={column name=$\|e_h\|_h$},
   columns/2/.style ={column name=$\operatorname{EOC}$},
        columns/3/.style ={column name=$|e_h|_{H^1(\Omega)}$},
   columns/4/.style ={column name=$\operatorname{EOC}$},
   columns/5/.style ={column name=$\|e_h\|_{L^2(\Omega)}$},
   columns/6/.style ={column name=$\operatorname{EOC}$},
         columns/7/.style ={column name=$\eta_h$},
   columns/8/.style ={column name=$\operatorname{EOC}$},
%   columns/10/.style ={column name=},
   every head row/.style={before row=\toprule,after row=\midrule},
every last row/.style={after row=\bottomrule}
   ]{table-unif-ndof_MA_05.csv}}
\end{subtable}
\bigskip
%%% s = 0.2, adapt
\caption{$p=4$, $a=0.5$, adaptive refinement}\label{table21}
\begin{subtable}{1\linewidth}
%\centering
 \resizebox{\linewidth}{!}{
\begin{filecontents*}{table-adapt_MA_05.csv}
8.100000000000000000e+01 8.885535614132523009e-04 0.000000000000000000e+00 3.135685197569323756e-05 0.000000000000000000e+00 2.558004085498390167e-06 0.000000000000000000e+00 9.901844463934765094e-04 0.000000000000000000e+00
1.390000000000000000e+02 6.299753498337940525e-04 -6.368489948760962527e-01 2.180056222736916787e-05 -6.731117937667596696e-01 1.734185624676687846e-06 -7.197621168331947983e-01 7.345876508230179611e-04 -5.529040975152208492e-01
2.750000000000000000e+02 2.815729572393676941e-04 -1.180261543044973660e+00 7.504057976534663177e-06 -1.563089922374386154e+00 2.955186884910400283e-07 -2.593534335336794783e+00 2.855803488230811648e-04 -1.384713216096425947e+00
4.950000000000000000e+02 9.875831022341262421e-05 -1.782476696565880259e+00 1.534328886811276547e-06 -2.700556093553666770e+00 5.567415427438074139e-08 -2.839833149536580503e+00 1.027647420615576083e-04 -1.738864054512081037e+00
1.107000000000000000e+03 3.640136270247941870e-05 -1.240066959714813732e+00 4.149587639966541879e-07 -1.624734179701296855e+00 1.302766103891489594e-08 -1.804608364039897328e+00 3.794388409385094285e-05 -1.237910835187092662e+00
2.531000000000000000e+03 1.107421892153654468e-05 -1.438987656482488831e+00 1.021574465872172335e-07 -1.694958076717399953e+00 2.092043731579948073e-09 -2.211632470897345382e+00 1.241950601116615753e-05 -1.350535592807786678e+00
4.915000000000000000e+03 4.528134420262304745e-06 -1.347506985608164287e+00 3.541085049510745720e-08 -1.596403820750745606e+00 6.215590831594879096e-10 -1.828698689026778901e+00 4.419234592102301345e-06 -1.556934121006759320e+00
9.285000000000000000e+03 1.478980076700783905e-06 -1.759066916241398770e+00 7.004726565545279079e-09 -2.547416500905537617e+00 7.541098112267987277e-11 -3.315909572296974694e+00 1.597547036884828302e-06 -1.599565569453346159e+00
\end{filecontents*}
\pgfplotstabletypeset[dec sep align,
   fixed zerofill,
   precision=3,
   col sep=space,
%      columns/0/.style ={column name=$h$},
   columns/0/.style ={int detect, column name=$\operatorname{ndofs}$},
   columns/1/.style ={column name=$\|e_h\|_h$},
   columns/2/.style ={column name=$\operatorname{EOC}$},
        columns/3/.style ={column name=$|e_h|_{H^1(\Omega)}$},
   columns/4/.style ={column name=$\operatorname{EOC}$},
   columns/5/.style ={column name=$\|e_h\|_{L^2(\Omega)}$},
   columns/6/.style ={column name=$\operatorname{EOC}$},
         columns/7/.style ={column name=$\eta_h$},
   columns/8/.style ={column name=$\operatorname{EOC}$},
%   columns/10/.style ={column name=},
   every head row/.style={before row=\toprule,after row=\midrule},
every last row/.style={after row=\bottomrule}
   ]{table-adapt_MA_05.csv}}
    \end{subtable}
\end{table}

\begin{table}
\caption{$p=4$, uniform refinement}\label{table22}
\begin{subtable}{1\linewidth}
 \resizebox{\linewidth}{!}{
\begin{filecontents*}{table-unif-ndofinc_MA_unknown.csv}
4.225000000000000000e+03 1.043738686394366955e+00 0.000000000000000000e+00 1.213943566493072730e-02 0.000000000000000000e+00 4.899840110642998160e-04 0.000000000000000000e+00 4.508703215015644750e-01 0.000000000000000000e+00
1.664100000000000000e+04 1.090160067837036229e+00 3.174334931448469926e-02 6.871417032961854923e-03 -4.151357467678785662e-01 1.237449757641177121e-04 -1.003865987803788107e+00 4.360238575434919328e-01 -2.442484234315126651e-02
6.604900000000000000e+04 1.034947501046359353e+00 -3.770243257394137459e-02 3.303043651501553959e-03 -5.313829870293115265e-01 2.496963086386766372e-05 -1.161066892790884930e+00 4.238100557017767112e-01 -2.061010087153335862e-02
\end{filecontents*}
\pgfplotstabletypeset[dec sep align,
   fixed zerofill,
   precision=3,
   col sep=space,
%      columns/0/.style ={column name=$h$},
   columns/0/.style ={int detect, column name=$\operatorname{ndofs}$},
   columns/1/.style ={column name=$\|\theta_h\|_h$},
   columns/2/.style ={column name=$\operatorname{EOC}$},
        columns/3/.style ={column name=$|\theta_h|_{H^1(\Omega)}$},
   columns/4/.style ={column name=$\operatorname{EOC}$},
   columns/5/.style ={column name=$\|\theta_h\|_{L^2(\Omega)}$},
   columns/6/.style ={column name=$\operatorname{EOC}$},
         columns/7/.style ={column name=$\eta_h$},
   columns/8/.style ={column name=$\operatorname{EOC}$},
%   columns/10/.style ={column name=},
   every head row/.style={before row=\toprule,after row=\midrule},
every last row/.style={after row=\bottomrule}
   ]{table-unif-ndofinc_MA_unknown.csv}}
\end{subtable}
%%% s = 0.1, uniform
\bigskip
\caption{$p=4$, adaptive refinement}\label{table24}
\begin{subtable}{1\linewidth}
%\centering
 \resizebox{\linewidth}{!}{
\begin{filecontents*}{table-adapt_MA-inc_unknown.csv}
1.985000000000000000e+03 1.042029368630612352e+00 0.000000000000000000e+00 1.192714594011413075e-02 0.000000000000000000e+00 4.373449089219525062e-04 0.000000000000000000e+00 4.518104530828573329e-01 0.000000000000000000e+00
3.751000000000000000e+03 1.090871106265715040e+00 7.197701651063159733e-02 6.915200235381452863e-03 -8.565242074837270403e-01 1.353192413251997695e-04 -1.843304297805511416e+00 4.366403136283053543e-01 -5.366557852046364291e-02
7.063000000000000000e+03 1.037800061469112745e+00 -7.880795342157095873e-02 3.415064196285200974e-03 -1.114843181230515734e+00 4.693588535925082147e-05 -1.673158767668226687e+00 4.241184819614405299e-01 -4.597778163141733737e-02
1.477500000000000000e+04 1.003066343400422555e+00 -4.612252073611472492e-02 1.808046486535515220e-03 -8.616421000082170600e-01 2.396810144888799829e-05 -9.105663672825300425e-01 4.131613050038088208e-01 -3.546396283909993008e-02
2.900900000000000000e+04 9.749677100746193537e-01 -4.211332436666966744e-02 8.518064966614217158e-04 -1.115572583408316110e+00 5.543659144011379877e-06 -2.170053521676039487e+00 4.034020877976972619e-01 -3.543109992472079361e-02
5.918900000000000000e+04 9.533533956952976451e-01 -3.143705382769906853e-02 4.245106029244425359e-04 -9.765720201621387320e-01 3.558010163423730545e-06 -6.218411636410646759e-01 3.944118991869394741e-01 -3.160436586807339326e-02
\end{filecontents*}
\pgfplotstabletypeset[dec sep align,
   fixed zerofill,
   precision=3,
   col sep=space,
%      columns/0/.style ={column name=$h$},
   columns/0/.style ={int detect, column name=$\operatorname{ndofs}$},
   columns/1/.style ={column name=$\|\theta_h\|_h$},
   columns/2/.style ={column name=$\operatorname{EOC}$},
        columns/3/.style ={column name=$|\theta_h|_{H^1(\Omega)}$},
   columns/4/.style ={column name=$\operatorname{EOC}$},
   columns/5/.style ={column name=$\|\theta_h\|_{L^2(\Omega)}$},
   columns/6/.style ={column name=$\operatorname{EOC}$},
         columns/7/.style ={column name=$\eta_h$},
   columns/8/.style ={column name=$\operatorname{EOC}$},
%   columns/10/.style ={column name=},
   every head row/.style={before row=\toprule,after row=\midrule},
every last row/.style={after row=\bottomrule}
   ]{table-adapt_MA-inc_unknown.csv}}
   \end{subtable}
\end{table}

%, and % Furthermore, we replace the error estimator with 
%$$\eta_h=\|\max_{\alpha\in\Lambda}\{\gamma^\alpha(A:D^2_hu_h-f^\alpha)\}\|_{L^2(\X)}+\|\tilde{g}-g_h\|_{L^2(\X)}+\left(\sum_{e\in\Ei}\frac{1}{h_e}\int_e\jump{\partial u_h/\partial n}_e^2\,ds\right)^\frac{1}{2},$$
%where $g_h\in V_h$ is the numerical approximation of the extended boundary data $\tilde g\in H^2(\X)$.
\begin{figure}[h]
\begin{center}
 \hspace{-1.2cm}%
 \scalebox{0.66}{
    \input{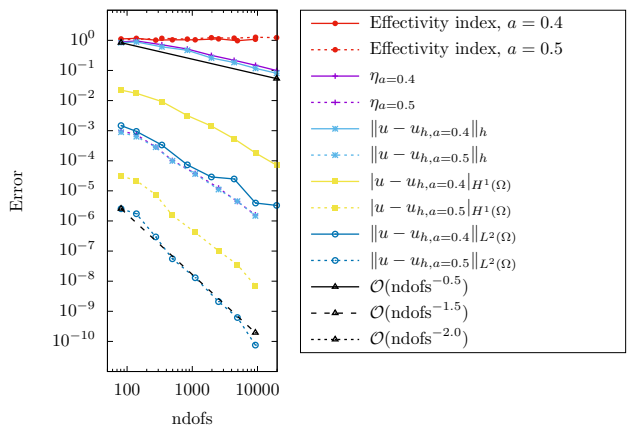}}
%    \unskip
%    }% P1 norec
%    &
%    \hspace{-0.7cm}
   %\vspace{20cm}
%    \scalebox{0.66}{%
 %   \renewenvironment{table}[1][]{\ignorespaces}{\unskip}%
 %  \input{exp_3-dualcons.tex}
%    \unskip
  %  }% P1 norec
  \end{center}
 \caption{Convergence rates for Experiment~\ref{exp4}, we observe faster convergence rates for $u_{0.5}$ than for $u_{0.4}$.}
\label{plot:exp4:1}
\end{figure}
%\subsection{Experiment 4}\label{exp5}
%In this experiment, we consider the following simply supported plate problem
%\begin{equation}\label{exp:5}
%\left\{
%\begin{aligned}
%\Delta^2\Phi&= \phi,\quad\mbox{in}\quad\X,\\
%\Delta\Phi=\Phi & = 0,\quad\mbox{on}\quad\partial\X,
%\end{aligned}
%\right.
%\end{equation}
%where $\X=(0,1)^2$. For this example, and the function $\phi$ is chosen so that the true solution $\Phi = \sin(\pi x)\sin(\pi y)$. We apply FEM~(\ref{FEMnew2}), with uniform mesh refinement, and use polynomial degrees $p=2,3,4$, and we consider the convergence of the method in the $\|\cdot\|_{h,\Delta}$-norm. We observe optimal convergence rates of order $\|\Phi-\Phi_h\|_{h,\Delta} = \mathcal{O}(h^{p-1})$, plotting the results in Figure~\ref{plot:exp5:1}. 
%\begin{figure}[h]
%\begin{center}
%%  \begin{tabular}{lrr}
%%  \hspace{-3cm}
%    \scalebox{0.7}{%
%%    \renewenvironment{table}[1][]{\ignorespaces}{\unskip}%
%    \input{exp_4.tex}
%%    \unskip
%    }% P1 norec
%%    &
%%    \hspace{0.1cm}
%%   %\vspace{20cm}
%%    \scalebox{0.265}{%
%%    \renewenvironment{table}[1][]{\ignorespaces}{\unskip}%
%%   \input{exp2refmesh.tex}
%%    \unskip
%%    }% P1 norec
%%  \end{tabular}
%  \end{center}
% \caption{Convergence rates for Experiment~\ref{exp5}.}
%\label{plot:exp5:1}
%\end{figure}
\end{section}
\begin{section}{Concluding Remarks}\label{sec:6}
In this paper, we were successful in proposing and analysing a $C^0$-interior penalty method for the approximation of the fully nonlinear Hamilton--Jacobi--Bellman equation with inhomogeneous Dirichlet boundary data. The analysis consisted of three parts: proving a stability estimate, a quasi-optimal a priori error estimate, and an a posteriori error estimate in a $H^2$-style norm. All of the aforementioned analysis was undertaken, assuming regularity no higher than $H^2(\X)$, as implied by the corresponding PDE theory. All of the theoretical results were confirmed in the experiments section, which included the implementation of an adaptive method, based upon the proven a posteriori error estimate. Furthermore, we were able to apply the proposed method to the fully nonlinear Monge--Amp\`{e}re equation, providing a uniquely solvable, optimally convergent, and adaptive finite element method.

%Alongside the first method, we proposed a second method that is adjoint consistent with a fourth-order anisotropic simply supported plate problem. The consistency holds in the case of additional regularity of the coefficients, and of the corresponding dual problem. For this method, we were successful in undertaking the same analysis as for the first method, as well as proving higher order convergence rates in the $L^2$-norm and $H^1$-seminorm, with rates determined by the index of elliptic regularity of the dual problem. A byproduct of this second method, and the analysis undertaken, is the proposition of a third method, which approximates solutions to the fourth-order dual problem; in this case we proved that the method is stable, and proved a quasi-optimal error estimate in a second norm equivalent to the $H^2$-style norm used for the aforementioned analysis for the first two methods. In the experiments section we confirmed the theoretical results for the adjoint consistent method, and the dual method.
 \end{section}
 %\newpage
%
%\input{table-exp3.tex}
%
\bibliographystyle{plain}
\bibliography{toskaweckiBIB}

Susanne C. Brenner, Department of Mathematics and Center for Computation and Technology, Louisiana State University, Baton Rouge, LA 70803, USA

\emph{E-mail address}: \texttt{brenner@math.lsu.edu}
\\

Ellya L. Kawecki, Department of Mathematics,
University College London, London, WC1H 0AY,
UK

\emph{E-mail address}: \tt{e.kawecki@ucl.ac.uk}
\end{document}